\documentclass[11pt, leqno]{amsart}
\usepackage{amsmath,amssymb,amsthm}
\usepackage[margin=1in]{geometry}
\usepackage[T1]{fontenc}
\usepackage{lmodern}
\usepackage{microtype}
\usepackage{graphicx}
\usepackage[colorlinks=true,linkcolor=blue,citecolor=blue,urlcolor=blue]{hyperref}
\hypersetup{pdftitle={Fractional Bombieri--De Giorgi--Giusti graphs},
pdfauthor={Manuel del Pino and Juncheng Wei},
pdfsubject={Fractional minimal graphs and inner--outer gluing}}
\numberwithin{equation}{section}
\theoremstyle{plain}
\newtheorem{theorem}{Theorem}[section]
\newtheorem{lemma}[theorem]{Lemma}
\newtheorem{proposition}[theorem]{Proposition}
\newtheorem{corollary}[theorem]{Corollary}
\theoremstyle{definition}
\newtheorem{definition}[theorem]{Definition}
\theoremstyle{remark}

\allowdisplaybreaks[4]

\newcommand{\R}{\mathbb R}
\newcommand{\eps}{\varepsilon}
\newcommand{\PV}{\operatorname{PV}}
\newcommand{\J}{\mathcal J}
\newcommand{\norm}[1]{\left\lVert#1\right\rVert}
\newcommand{\abs}[1]{\left\lvert#1\right\rvert}
\newcommand{\dd}{\,\mathrm d}

\newcommand{\sgn}{\operatorname{sgn}}

\newcommand{\HH}{\mathcal H}

\newcommand{\cP}{\mathcal P}
\allowdisplaybreaks[2]
\title[Fractional Bombieri--De Giorgi--Giusti graphs]{Fractional Bombieri--De Giorgi--Giusti graphs}
\author[M. del Pino]{Manuel del Pino}
\address{Manuel del Pino, Department of Mathematical Sciences, University of Bath, Bath BA2 7AY, United Kingdom}
\email{mdp59@bath.ac.uk}

\author[J.C. Wei]{Juncheng Wei}
\address{Juncheng Wei, Department of Mathematics, Chinese University of Hong Kong, Shatin, New Territories, Hong Kong}
\email{wei@math.cuhk.edu.hk}

\begin{document}
\begin{abstract}
  We construct smooth nonaffine entire graphs in $\R^9$ with zero fractional
  mean curvature, for every perimeter exponent sufficiently close to one.
  The construction is an inner--outer gluing. Exact fractional minimal
  hypersurfaces asymptotic to the Simons cone provide the outer foliation;
  two explicit Jacobi fields on the cylindrical cone determine its matching
  with the inner region.  In a slowly expanding compact set we insert a compatible classical Bombieri--De Giorgi--Giusti graph. A positive global barrier then gives a
  uniform weighted inverse and a small normal deformation yields the
  required graphs. Their growth is asymptotically homogeneous of degree
  converging to three, and they converge locally to a classical
  Bombieri--De Giorgi--Giusti graph. As far as we know this seems to be the first {\em nontrivial} fractional  minimal graph.
\end{abstract}
\maketitle

\medskip
\noindent\textbf{Mathematics Subject Classification (2020).}
53A10, 49Q05, 35R11.

\noindent\textbf{Keywords.} Fractional minimal surfaces, Bernstein problem,
Simons cone, inner--outer gluing.

\section{Introduction}
\label{sec:introduction}

Let $ 0<\sigma<1$. In a seminal  paper \cite{CRS2010}, Caffarelli, Roquejoffre and Savin introduced the $\sigma$-fractional perimeter whose critical points are fractional minimal surfaces satisfying the zero fractional mean curvature equation. More precisely let $E$ be an open subset of $\R^n$. For
$0<\sigma<1$, the fractional mean curvature of $\partial E$ at a point
$X\in\partial E$ is
\begin{equation}\label{eq:physicalcurvature}
  H_E^\sigma(X)
  =\PV\int_{\R^n}
  \frac{\chi_{E^c}(Y)-\chi_E(Y)}{|X-Y|^{n+\sigma}}\dd Y.
\end{equation}

The regularity of fractional minimal surfaces is a challenging problem (see \cite{BFV, CV2013, DSV2025, FV2017} and references therein). Asymptotically, as \(\sigma \to 1^-\), the fractional perimeter recovers the classical perimeter, allowing many classical minimal surface results to extend to the fractional setting \cite{CCS2020, CV2013, NP, Pa}. Conversely, as \(s \to 0^+\), the fractional perimeter approaches the relative volume \cite{Ca2025, DFPV, DSTV, Tsi}. For further results on stable fractional minimal surfaces, we refer to \cite{DDPW, Ca2020, Ca2025, CSSV2026, DSTV, SV} and references therein.

The main purpose of this paper is to extend the classical Bombieri-De Giorgi-Giusti graph \cite{BDG} to the fractional setting, as $ \sigma$ is close to 1. Let $E$ be an open subset of $\R^9$ with smooth boundary. We are interested in boundaries of the form
\[
\Gamma_F=\{(x,F(x)):x\in\R^8\},
\qquad E_F=\{(x,t):t<F(x)\},
\]
for which the fractional mean curvature $H_{E_F}^\sigma=0$.  Trivial solutions are affine functions, and the question addressed in this paper is whether
there are non-affine entire solutions.

Let us recall that for $\sigma=1$, Bombieri, De Giorgi and Giusti \cite{BDG} constructed a nonplanar
entire minimal graph over $\R^8$. The geometry of their graph is
closely tied to the Simons cone
\begin{equation*}
  C=\{x=(u,v)\in\R^4\times\R^4:|u|=|v|\}.
\end{equation*}
Write $|u|=r\cos\theta$, $|v|=r\sin\theta$, with
$\pi/4<\theta<\pi/2$ on the side $|v|>|u|$. The classical homogeneous model is
\begin{equation}\label{eq:classicalprofile}
  t=F_0(r,\theta)=r^3g(\theta),
\end{equation}
where $g(\pi/4)=0$, $g'(\pi/2)=0$, $g>0$ in the sector and
$g'(\pi/4)=1$; see \cite[Lemma 2.1 and Theorem 2]{DPKW}.
The equation for this model is \cite[(2.6)]{DPKW}.
Far from the origin, the horizontal sections of the graph are close to
dilations of a minimal hypersurface asymptotic to $C$, whereas near
$C\times\R$ the graph is better described as a small normal
displacement of this cylindrical cone. Both descriptions will be
needed in the fractional problem.

There are two difficulties  in the perturbation argument.  The
first concerns the behavior at infinity. On every fixed compact set, a
classical minimal graph has small normalized fractional mean curvature
when $\sigma$ is close to one. This does not, however, make the
curvature small in the weighted space in which the equation has to be
solved on the whole graph: the leading error at infinity has the slow
spatial decay, and expanding to higher order in $1-\sigma$ improves its
coefficient but not its decay. The exterior geometry must therefore be
constructed for the exact value of $\sigma$. This is also the difficulty
in the construction of the fractional catenoid in
\cite[Section 2, (2.3)--(2.6)]{DDPW}.

The second difficulty comes from the cylindrical cone. The radial
Jacobi equation on $C$ has the two radial powers $r^{-m_\sigma}$ and
$r^{-n_\sigma}$, where
\begin{equation}\label{eq:indicialnotation}
  m_\sigma\longrightarrow2,
  \qquad n_\sigma\longrightarrow3,
  \qquad m_\sigma+n_\sigma=6-\sigma.
\end{equation}
The exact exterior leaf contains both powers in its expansion, and
retaining only the first of them in the interior leaves a matching
error which is not small enough. We therefore construct explicit
solutions of the cylindrical Jacobi equation which account both for the
second radial power and for the additional coefficient of $r^{-m_\sigma}$
in their expansion.

Throughout the paper $\sigma \in (0,1)$ denotes the perimeter exponent in
\eqref{eq:physicalcurvature}, and the surface Jacobi operator has order
$\tau$, where
\begin{equation*}
  \tau=1+\sigma,\qquad \eps=1-\sigma=2-\tau.
\end{equation*}
We fix $\beta=5/2$ throughout the paper.
Our main result is the following.

\begin{theorem}\label{thm:main}
  There exists $\sigma_0\in(0,1)$ such that, for every
  $\sigma\in(\sigma_0,1)$, there is a smooth nonaffine function
  $F_\sigma:\R^8\to\R$ whose graph has vanishing $\sigma$-mean curvature
  and which satisfies
  \[
  F_\sigma(Uu,Vv)=F_\sigma(u,v),\qquad
  F_\sigma(v,u)=-F_\sigma(u,v)
  \quad\hbox{for all }U,V\in O(4).
  \]
  Moreover, there are an exponent $p_\sigma>1$, with $p_\sigma\to3$ as
  $\sigma\uparrow1$, and a nonzero continuous function
  $F_{\infty,\sigma}:\R^8\to\R$, homogeneous of degree $p_\sigma$, such
  that
  \begin{equation}\label{eq:mainasymptotic}
    R^{-p_\sigma}F_\sigma(Rx)\longrightarrow F_{\infty,\sigma}(x)
    \quad\hbox{as }R\longrightarrow\infty,
  \end{equation}
  locally uniformly on $\R^8\setminus C$. Finally, the functions
  $F_\sigma$ can be chosen so that, as $\sigma\uparrow1$, they converge
  locally in $C^2$ to a classical entire minimal graph $F_*$ with cubic
  growth.
\end{theorem}

The limiting profile $F_{\infty,\sigma}$ is determined by the family of
dilations of a fractional minimal hypersurface asymptotic to $C$, and
the construction gives $p_\sigma=1+m_\sigma$, where $m_\sigma$ is the
slowly decaying radial indicial root in \eqref{eq:indicialnotation}.

In the same coordinates the fractional limiting profile is
\begin{equation}\label{eq:fractionalprofile}
  F_{\infty,\sigma}(r,\theta)=r^{p_\sigma}g_\sigma(\theta).
\end{equation}
The angular function is normalized by
\begin{equation*}
  \begin{gathered}
    g_\sigma(\pi/4)=0,\qquad g_\sigma'(\pi/4)=1,\qquad
    g_\sigma'(\pi/2)=0,\\
    g_\sigma>0\quad\hbox{on }(\pi/4,\pi/2),\qquad
    g_\sigma(\pi/2-\theta)=-g_\sigma(\theta).
  \end{gathered}
\end{equation*}
We give its equation explicitly. Fix $e=(1,0,0,0)\in S^3$ and write
\[
x_\theta=(\cos\theta\,e,\sin\theta\,e),\qquad
y=\rho(\cos\varphi\,\omega,\sin\varphi\,\zeta),
\quad \omega,\zeta\in S^3.
\]
Set
\[
D^2=|x_\theta-y|^2
=1+\rho^2-2\rho(\cos\theta\cos\varphi\,\omega_1
+\sin\theta\sin\varphi\,\zeta_1).
\]
Then $g_\sigma$ satisfies the nonlocal angular equation
\begin{equation}\label{eq:angular-fractional}
  \lim_{\delta\downarrow0}\int_0^\infty\!\int_0^{\pi/2}
  \!\int_{S^3}\!\int_{S^3}
  \mathbf1_{\{D>\delta\}}
  \frac{\sgn\big(\rho^{p_\sigma}g_\sigma(\varphi)
  -g_\sigma(\theta)\big)}{D^{8+\sigma}}
  \rho^7\cos^3\varphi\sin^3\varphi
  \dd\omega\,\dd\zeta\,\dd\varphi\,\dd\rho=0.
\end{equation}
Indeed, $x_\theta$ belongs to the level set
$\{F_{\infty,\sigma}=g_\sigma(\theta)\}$, which is a dilation of
$M_\sigma$. With $E_a=\{F_{\infty,\sigma}<a\}$, the numerator is
$\chi_{E_a^c}-\chi_{E_a}$, and the measure in these polar coordinates
is $\rho^7\cos^3\varphi\sin^3\varphi\dd\rho\dd\varphi\dd\omega\dd\zeta$.
Thus \eqref{eq:angular-fractional} is exactly the fractional minimal
surface equation for the level sets. The exclusion $D>\delta$ is
made before any angular integration.

Multiply the curvature in \eqref{eq:angular-fractional} by
$-c_{7,\sigma}/2$, where $c_{7,\sigma}$ is given in
\eqref{eq:normalization}, and let $\sigma\uparrow1$. On compact
subintervals of $(\pi/4,\pi/2)$ the limit $g_1$ satisfies
\begin{equation}\label{eq:angular-classical}
  \left(\frac{g_1'(\theta)\sin^3(2\theta)}
  {\sqrt{9g_1(\theta)^2+g_1'(\theta)^2}}\right)'
  +\frac{21g_1(\theta)\sin^3(2\theta)}
  {\sqrt{9g_1(\theta)^2+g_1'(\theta)^2}}=0,
  \qquad \frac\pi4<\theta<\frac\pi2.
\end{equation}
This is the same differential equation as \cite[(2.6)]{DPKW}.
It follows from
$\operatorname{div}(DF_{\infty,1}/|DF_{\infty,1}|)=0$, since
$F_{\infty,1}=r^3g_1(\theta)$ has classical minimal level sets.
For $\sigma<1$ the equation for $g_\sigma$ is
\eqref{eq:angular-fractional}, with all angular values coupled by
the integral.

On a fixed level $F_{\infty,\sigma}=A$, the normalization
$g_\sigma'(\pi/4)=1$ gives
$\theta-\pi/4\sim Ar^{-p_\sigma}$ as $r\to\infty$. Its normal height
is therefore asymptotic to $Ar^{1-p_\sigma}$. Comparing this with
$Ar^{-m_\sigma}$ in the leaf expansion gives $p_\sigma=1+m_\sigma$.
Thus
\begin{equation*}
  I_\sigma(p_\sigma-1)=0,\qquad
  p_\sigma=1+m_\sigma=7-\sigma-n_\sigma,\qquad
  p_\sigma=3+o(1)\quad(\sigma\uparrow1),
\end{equation*}
where $I_\sigma$ is the explicit symbol in \eqref{mod:symbol} and
$m_\sigma$ is its root near two. In particular $n_\sigma$ in
\eqref{eq:indicialnotation} is a decay exponent, whereas $p_\sigma$
is the degree of $F_{\infty,\sigma}$. The sign of $p_\sigma-3$ is
the sign of $I_\sigma(2)$ for $\sigma$ close to one. The convergence
of the roots proved here does not determine this sign.

The exterior approximation will be written as
$A_{\rm out}(t)=r^{p_\sigma}g_\sigma(\theta)$, with
$A_{\rm out}(t)=t+o(t)$. This relation determines the horizontal
sections. Near $C\times\R$ we use normal coordinates over the
cylinder and solve its Jacobi equation. The function $A_{\rm out}$
is chosen so that the two descriptions agree.

\subsection{Relation with the fractional catenoid}\label{intro:catenoid}

The construction follows \cite[Section 2]{DDPW}, with $s=\sigma$
and $\eps=1-\sigma$. For the catenoid the interaction of its two ends
leads to
\[
\Delta f_\eps=\frac{\eps}{f_\eps},
\]
and changes the classical logarithmic end into a conical end
\cite[(2.5)--(2.8)]{DDPW}. After constructing this approximation,
one estimates its error, solves the linearized equation and applies a
contraction argument. The counterparts of
\cite[Propositions 2.1--2.3]{DDPW} are, respectively,
Proposition~\ref{res:assembled} with Lemma~\ref{glu:diagonal},
the inverse estimate \eqref{glu:inverseestimate}, and
Proposition~\ref{glu:quadratic}.

Several calculations are used directly. The curvature variation is
\cite[Proposition B.1]{DDPW}; the sign and normalization are fixed
in Subsection~\ref{sec:preliminaries}. The radial power calculation on
$C$ is \cite[Section 11.3, (11.7)]{DDPW}, and its symmetry is
\cite[Remark 11.1]{DDPW}. The local nonlinear estimate uses
\cite[Lemmas 8.3--8.5]{DDPW}, with the H\"older exponents related as
explained below. These calculations are cited where they are applied.

The exterior geometry is different. The horizontal sections used here
are exact fractional minimal hypersurfaces. Their vertical dependence,
the second coefficient in the cylindrical expansion, and the cutoff
between the two approximations determine the error. These terms are
estimated in Section~\ref{res:section}. The catenoid inverse uses the
translation and dilation fields in \cite[Section 7]{DDPW}; here we use
positive barriers and the identity in the proof of
\cite[Corollary B.1]{DDPW}.

There is also a difference in the small-parameter estimates. In
\cite[(2.13), (2.15)]{DDPW} the normalized error is
$O(\sqrt\eps/|\log\eps|)$ and the nonlinear estimate has a factor
$\eps^{-1/2}$. Here the error tends to zero by the choice of the compact
matching radius, and the nonlinear constant is uniform. The geometric
reason and the required integrals are given in
Subsection~\ref{glu:nonlinearsection}.

\subsection{Outline of the construction}

We first construct a rotationally invariant fractional minimal
hypersurface $M_\sigma\subset\R^8$ on one side of $C$. Its normal
height has the expansion
\begin{equation*}
  h_\sigma(r)=a_0r^{-m_\sigma}+b_\sigma r^{-n_\sigma}+N_\sigma(r),
  \qquad N_\sigma(r)=O(r^{-6}),
\end{equation*}
with smooth bounds and fixed $a_0>0$. The starting classical
smoothing and its expansion are given by
\cite[Theorem~2.1]{HS} and \cite[Proposition~3.3]{Szekelyhidi}.
We parametrize the signed dilations by
\[
A=a_0\sgn(\lambda)|\lambda|^{p_\sigma},\qquad
p_\sigma=1+m_\sigma.
\]
Their first two height terms are $Ar^{-m_\sigma}$ and
$b\phi_q(A)r^{-n_\sigma}$, where
$q=(1+n_\sigma)/p_\sigma$ and $b=b_\sigma a_0^{-q}$.

To match these leaves near $C\times\R$, we construct cylindrical
Jacobi fields from the Poisson kernel of
Section~\ref{mod:cylindersection}. For $d=n_\sigma-m_\sigma$, the
field associated with $\phi_q(t)=\sgn(t)|t|^q$ has the two terms
\[
r^{-n_\sigma}\phi_q(t)
+k(q,d)r^{-m_\sigma}\phi_{q-d}(t).
\]
The second term determines a correction to the exterior parameter
$A$. Repeating the matching once gives two fields, which are joined
to the exact leaves on $|t|\asymp r^{5/2}$. Cutting off $N_\sigma$
near $t=0$ gives a smooth approximate graph $\cP_\sigma$.

At $\sigma=1$, Dirichlet problems with these boundary values select
an entire classical graph $\Gamma_*$. Comparison uses the BDG graph
of \cite{BDG}, in the form of \cite[Theorem~2]{DPKW}, and strict
exterior barriers, and gives a normal
difference from $\cP_1$ of weight $\beta_1=31/12>\beta$.
We cut off this remainder at radius $R$ and insert it into
$\cP_\sigma$. The completed approximation satisfies
\begin{equation*}
  \norm{\HH_\sigma(\Sigma_{\sigma,R})}_{Y_{\sigma,\beta}}
  \le CR^{-1/12}+\omega_R(1-\sigma),\qquad
  \omega_R(\eps)\longrightarrow0
  \quad\hbox{for fixed }R.
\end{equation*}
Choose first $R$ large and then $\sigma$ close to one; a diagonal
choice gives $R(\sigma)\to\infty$ and vanishing error. The cutoff
at $r\asymp R$ is distinct from the matching on
$|t|\asymp r^{5/2}$, which is performed at the exact fractional
parameter on the unbounded end.

A global positive barrier gives a uniformly bounded Jacobi inverse.
The quadratic curvature estimate then yields a small normal correction
by contraction. Its horizontal projection is a global diffeomorphism,
so the resulting surface is an entire graph.

\subsection{Organization of the paper}

Subsection~\ref{sec:preliminaries} fixes notation and weighted spaces.
Sections~\ref{mod:radial} and~\ref{mod:cylindersection} construct the
fractional leaf and cylindrical fields. Section~\ref{res:section}
estimates the curvature of the matched graph, and
Section~\ref{glu:coresection} inserts the classical core.
Section~\ref{glu:inverse} proves the inverse and nonlinear estimates
and completes the proof of Theorem~\ref{thm:main}.

\subsection{Notation and analytic conventions}\label{sec:preliminaries}

Put $E_C=\{|u|>|v|\}$. For $P=(u,v)\in C$, $r=|P|$, its outward normal is
\begin{equation*}
  \nu_C(P)=r^{-1}(-u,v).
\end{equation*}
The link is $S^3(1/\sqrt2)\times S^3(1/\sqrt2)$ and the radial measure
is $r^6\dd r$. In the coordinates $x=P+s\nu_C(P)$,
\begin{equation*}|A_C|^2=6r^{-2},\quad
  J_C^1=\partial_r^2+6r^{-1}\partial_r+6r^{-2},\quad
  J_C(P,s)=(1-s^2/r^2)^3.
\end{equation*}
The polar variables above satisfy
\[
|x|^2=r^2+s^2,\qquad \tan(\theta-\pi/4)=s/r.
\]
Thus $r$ and $|x|$ are comparable when $|s|\le cr$.
Thus $J_C^1r^{-z}=(z-2)(z-3)r^{-z-2}$. The set $E_C$ is $s<0$.
Moving its boundary to $s=h$ gives first variation $-h\delta_C$ for
$\chi_{E^c}$, and $-2h\delta_C$ for $\chi_{E^c}-\chi_E$.

For graphs we use the upward normal
\begin{equation*}
  \nu_F=\frac{(-DF,1)}{\sqrt{1+|DF|^2}}.
\end{equation*}
The set $\{s<h(P,t)\}$ describes the region above the graph
when $h_t>0$, and hence has the opposite orientation. Replacing $E$
by $E^c$ changes the sign of curvature, without changing its zero set.
The Jacobi operator is unchanged by simultaneous normal reversal.

For $M=\partial E\subset\R^N$, the first variation of
\eqref{eq:physicalcurvature} under $X\mapsto X+h\nu_M$ is
$-2\J_M^\sigma h$. This is \cite[Proposition B.1, (B.1)]{DDPW}
with the sign in \eqref{eq:physicalcurvature}. The operator is \cite[(1.5)]{DDPW}
\begin{equation*}
  \J_M^\sigma h(X)=\PV\int_M
  \frac{h(Y)-h(X)+(1-\nu_M(X)\cdot\nu_M(Y))h(X)}{|X-Y|^{N+\sigma}}
  \dd S_Y.
\end{equation*}
We fix the normalizations
\begin{equation}\label{eq:normalization}
  c_{7,\sigma}=\frac{14(1-\sigma)}{|S^6|},\quad
  D_\sigma=\int_\R(1+z^2)^{-(9+\sigma)/2}\dd z
  =\frac{\sqrt\pi\Gamma((8+\sigma)/2)}{\Gamma((9+\sigma)/2)},\quad
  \kappa_\sigma=\frac{c_{7,\sigma}}{D_\sigma}.
\end{equation}
On surfaces in $\R^9$ put
\begin{equation*}
  \HH_\sigma=-\frac{\kappa_\sigma}{2}H_E^\sigma,
  \qquad L_{M,\sigma}=\kappa_\sigma\J_M^\sigma;
\end{equation*}
in $\R^8$ replace $\kappa_\sigma$ by $c_{7,\sigma}$.
To make the comparison with \cite{DDPW} explicit, their parameter is
$s=\sigma$ and their curvature is
\begin{equation}\label{eq:DDPW-conventions}
  H_{\rm DDPW}^\sigma
  =\PV\int_{\R^N}\frac{\chi_E(Y)-\chi_{E^c}(Y)}{|X-Y|^{N+\sigma}}\dd Y
  =-H_E^\sigma,\qquad
  \HH_\sigma=\frac{\kappa_\sigma}{2}H_{\rm DDPW}^\sigma.
\end{equation}
Hence their moving-point expansion \cite[(2.9), (B.1)]{DDPW} becomes
$\HH_\sigma(M[h])=E+L_{M,\sigma}h+Q_\sigma(h)$, with
$L_{M,\sigma}=\kappa_\sigma\J_M^\sigma$ and
$Q_\sigma=(\kappa_\sigma/2)N$, where $N$ is their remainder.
The same conversion uses $c_{7,\sigma}$ for the leaf in $\R^8$.

We give an explicit formula for this remainder. Let $M=\partial E\subset
\R^9$, with outward normal $\nu$, and set $A_M=D_M\nu$. Assume its
normal deformation is embedded and $|uA_M|<1/2$. Write
\[
X_u=X+u(X)\nu(X),\qquad Z=Y-X,\qquad
W=u(Y)\nu(Y)-u(X)\nu(X).
\]
The outward normal of $M[u]$, multiplied by the area Jacobian of
$Y\mapsto Y+u(Y)\nu(Y)$, is
\begin{equation*}
  \begin{aligned}
    n_u&=\det(I+uA_M)
    \big[\nu-(I+uA_M)^{-1}\nabla_Mu\big],\\
    n_1&=u\,\operatorname{tr}(A_M)\nu-\nabla_Mu.
  \end{aligned}
\end{equation*}
The determinant and inverse act on the eight-dimensional tangent
space; $n_1$ is the term linear in $u$. The identity follows by taking
the normal to $(I+uA_M)v+(\nabla_Mu\cdot v)\nu$ and multiplying by
the area Jacobian.
The divergence of $(Y-X)|Y-X|^{-9-\sigma}$ is
$-\sigma|Y-X|^{-9-\sigma}$. Applying the divergence theorem to the
$E$ and $E^c$ outside a small ball gives
$H_E^\sigma(X)=(2/\sigma)\int_M Z\cdot\nu(Y)|Z|^{-9-\sigma}\dd S_Y$:
the leading half-sphere fluxes cancel, and their error is
$O(\delta^{1-\sigma})$. The flux at infinity tends to zero. Hence
\begin{equation*}
  \HH_\sigma(M[u])(X_u)
  =-\frac{\kappa_\sigma}{\sigma}\int_M
  \frac{(Z+W)\cdot n_u(Y)}{|Z+W|^{9+\sigma}}\dd S_Y.
\end{equation*}
Subtracting its constant and linear terms gives the explicit expression
\begin{equation}\label{eq:explicit-Q}
  \begin{aligned}
    Q_\sigma(u)(X)=-\frac{\kappa_\sigma}{\sigma}\int_M\bigg[
    &\frac{(Z+W)\cdot n_u(Y)}{|Z+W|^{9+\sigma}}
    -\frac{Z\cdot\nu(Y)+W\cdot\nu(Y)+Z\cdot n_1(Y)}{|Z|^{9+\sigma}}\\
    &+(9+\sigma)\frac{(Z\cdot\nu(Y))(Z\cdot W)}{|Z|^{11+\sigma}}
    \bigg]\dd S_Y.
  \end{aligned}
\end{equation}
All terms use the same principal-value truncation. For $C^2$ data the
numerators of the curvature term and of its combined linear terms
vanish to second order at $Y=X$, while
$(Z\cdot\nu(Y))(Z\cdot W)=O(|Z|^4)$. The local radial majorant is
$C\rho^{-\sigma}\dd\rho$. Thus these combined integrals are
absolutely convergent for $\sigma<1$. The first derivative is
$L_{M,\sigma}u$ by \cite[Proposition B.1]{DDPW}, so
\eqref{eq:explicit-Q} is precisely the remainder in the curvature
expansion above. The weighted estimates below control its
integrals at infinity. For the leaf in $\R^8$ replace $9$ by $8$,
$11$ by $10$, and $\kappa_\sigma$ by $c_{7,\sigma}$; the tangent
determinant then has dimension seven.

We write $\widehat\J_C^\sigma=c_{7,\sigma}\J_C^\sigma$ and
$\widehat\J_{C\times\R}^\sigma=\kappa_\sigma\J_{C\times\R}^\sigma$.
Integration in the vertical variable gives the exact identity
\begin{equation*}
  \widehat\J_{C\times\R}^\sigma f(P,t)=\widehat\J_C^\sigma f(P)
  \quad\hbox{if }f\hbox{ is independent of }t.
\end{equation*}
These normalized operators converge locally to $\Delta_M+|A_M|^2$;
the curvature convention is $\HH_1=-\operatorname{div}_M\nu_M$.

On a completed graph set $\mathfrak r(X)=(1+|x|^2)^{1/2}$.
Let $\Phi_X:B_2\to M$ be a tangent chart at scale $c\mathfrak r(X)$,
with uniform scaled smooth bounds. The auxiliary weights are
\begin{equation*}
  b_s=\frac94,\qquad\beta_1=\frac{31}{12},\qquad
  \gamma=\frac{21}{8},\qquad\eta=\frac18,
\end{equation*}
and fix a sufficiently small $\alpha>0$, with $\alpha<1/48$.
Restricting $\sigma$ close to one, we have
\begin{equation*}
  \tau+\alpha=2+\alpha_\sigma,\quad
  \alpha/2<\alpha_\sigma<\alpha<\eta,\qquad
  m_\sigma<b_s<\beta<\beta_1<\gamma<n_\sigma.
\end{equation*}
The weighted norms are
\begin{align*}
  \|h\|_{X_{\sigma,b}(M)}
  &=\sup_X\mathfrak r(X)^b
  \|h\circ\Phi_X\|_{C^{\tau+\alpha}(B_1)},\\
  \|f\|_{Y_{\sigma,b}(M)}
  &=\sup_X\mathfrak r(X)^{b+\tau}
  \|f\circ\Phi_X\|_{C^\alpha(B_1)}.
\end{align*}
At $\sigma=1$ use $C^{2,\alpha}$ in the first norm. The leaf norms use
seven-dimensional charts. In particular, $\|h\|_{X_{\sigma,b}}\le C$
implies $|\nabla^jh|\le C\mathfrak r^{-b-j}$ for $j\le2$, with the
scaled $C^{\alpha_\sigma}$ bound for $\nabla^2h$.
The notation $O_{C^k_{\rm sc}}(r^{-b})$ includes scaled derivatives up
to order $k$. On ends we usually write $r$ in place of $\mathfrak r$.

These powers agree with the geometric scaling. Curvature has degree
$-\sigma$, whereas its derivative with respect to a normal height has
degree $-1-\sigma=-\tau$. Rescaling a chart by $r$ therefore replaces
the physical height $h$ by $h/r$, the curvature by $r^\sigma\HH_\sigma$,
and the Jacobi operator by $r^\tau L_{M,\sigma}$. This explains the
difference of $\tau$ between the two norm weights. Throughout the
argument the correction is measured as a normal height.

We record the common local estimates used in the two perturbations.
The classical limits are \cite[Appendix A, Lemmas A.1, A.2 and A.4]{DDPW}.
For the surface integrals the hypotheses include bounded $C^{2,\alpha}$
test functions and a uniformly bounded remote integral
$\int_{M\setminus B_1(X)}|X-Y|^{-d-1-\sigma}\dd S_Y$, where
$d=\dim M$. We apply these results after localization. The radial
and area estimates below supply the weighted tails, also when the
test function is a power on an end. Thus the cited local limits are
used with their integrability hypotheses checked in each geometry.

The graph calculation is the one in \cite[Section 8]{DDPW}:
the implicit graph map (8.1), its first two variations in Lemma 8.3,
and the chain rule preceding that lemma. After integrating in the
normal variable, its local graph formula is, up to normalization
and orientation,
\[
\int\frac{\zeta(z)}{|z|^{d+\sigma}}
\left[\psi_{d,\sigma}\left(\frac{g(y+z)-g(y)}{|z|}\right)
-\psi_{d,\sigma}\left(\frac{Dg(y)\cdot z}{|z|}\right)\right]\dd z,
\qquad
\psi_{d,\sigma}(a)=\int_0^a(1+t^2)^{-(d+1+\sigma)/2}\dd t.
\]
This is the secant-slope identity on pp.~42--43 of \cite{DDPW}; its
dimension-dependent exponent follows from the ambient kernel.
The chart slopes and their scaled derivatives are bounded here, and
the implicit map has uniformly bounded first two variations on a
fixed small height ball. The elementary differentiations in
\cite[(8.3)--(8.6), Lemmas 8.3--8.4]{DDPW} therefore give the local
Taylor structure in dimensions seven and eight. Their numerical
weighted estimates concern the catenoid. We verify below the local
estimate from $C^{\tau+\alpha}$ to $C^\alpha$ needed here.

In a smooth scaled chart, the displacement kernel has the form
\[
(2-\tau)\frac{j_G(y+z)}{|G(y+z)-G(y)|^{d+\tau}},
\qquad d=7\ \hbox{or}\ 8.
\]
Its even part is uniformly elliptic; its first two displacement
derivatives have bounds $C(2-\tau)|z|^{-d-\tau-j}$.
The odd part gains one displacement power. Its operator is bounded
from $C^{1,\alpha}$ to $C^\alpha$, uniformly in $\tau$: pair the
increments and integrate $(2-\tau)|z|^{2-d-\tau}$ locally.
Interpolation bounds this term by
$\delta\|h\|_{C^{\tau+\alpha}}+C_\delta\|h\|_\infty$,
and the first term is absorbed in the local estimate.
The even amplitude is $C^\alpha$ in the center, uniformly after
scaling. The local Jacobi potential has the same regularity because
$1-\nu(X)\cdot\nu(Y)=O(|X-Y|^2/r^2)$.
In particular, the center differences satisfy the integral condition
\cite[(1.7)]{JX}, since
\[
(2-\tau)\int_{\R^d}|z|^{-d-\tau}
\min(|z|^2,\rho^2)\dd z\le C\rho^{2-\tau}.
\]
Let $\bar\alpha$ be the uniform exponent in that theorem. Choose
$\alpha<\bar\alpha/2$ and restrict $1-\sigma<\alpha/2$.
Then $\tau+\alpha-2\ge\alpha/2$ and
$\tau+\bar\alpha-2\ge\bar\alpha/2$.
The Schauder estimate \cite[Theorem 1.1]{JX}, followed
by localization, therefore yields
\begin{equation}\label{glu:schauder}
  \|h\|_{X_{\sigma,b}}
  \le C\big(\|L_{M,\sigma}h\|_{Y_{\sigma,b}}
  +\|h\|_{C_b^0}\big),\qquad
  \|h\|_{C_b^0}=\sup_M\mathfrak r^b|h|,
\end{equation}
whenever the weighted remote integrals are bounded at the same scale.
No unweighted bound on a rescaled copy of the whole surface is used.

The local remainder and its first two height variations are treated by
the second-difference calculation in
\cite[Lemmas 8.4--8.5 and the proof of Lemma 8.1]{DDPW}.
In those estimates take their H\"older exponent
$\alpha_{\rm DDPW}=\alpha_\sigma=\alpha-\eps$. Their input space
$C^{2,\alpha_{\rm DDPW}}$ is our $C^{\tau+\alpha}$, and their source
exponent $\alpha_{\rm DDPW}+\eps$ is our $\alpha$; see
\cite[(2.11)--(2.12), Lemma 8.5]{DDPW}.
Since $\alpha_\sigma\in[\alpha/2,\alpha]$, the exponents remain in a
fixed compact interval. The normalization factor, comparable to
$\eps$, cancels the $1/\eps$ in the local estimates of
\cite[Lemmas 8.1 and 8.4]{DDPW}; their radial integration is
$\eps\int_0^1s^{-\sigma}\dd s=1$.
Changing the dimension only changes $\psi_{d,\sigma}$ above: its
derivatives are uniformly bounded on the bounded slope range, and the
same Taylor expansion and estimates of differences at two evaluation
points apply.
This verifies the local estimate in our norms without importing the
catenoid's global weighted constant. Remote terms are differentiated
before integration and controlled by the separate absolute moments.

All unspecified constants are uniform for $\sigma$ close to one.
Dependence on a fixed compact gluing radius is indicated by $C_R$.

\section{The fractional smoothing of the Simons cone}\label{mod:radial}
\subsection{The radial equation}

We write $\tau=1+\sigma$, $\eps=1-\sigma$, and use the normalization
of Section~\ref{sec:preliminaries}. The cone $C=C_{4,4}$ is
$\sigma$-minimal for every $0<\sigma<1$. This is the case of two
four-dimensional factors in \cite[Theorem 3]{DDPW}; the interchange
$(u,v)\mapsto(v,u)$ gives opening parameter one. We use its Jacobi operator and radial-power
calculation from \cite[Section 11]{DDPW}. The normalized operator on
the seven-dimensional cone is
\begin{equation*}
  \widehat\J_C^\sigma=c_{7,\sigma}\J_C^\sigma,
  \qquad c_{7,\sigma}=\frac{14\eps}{|S^6|}.
\end{equation*}
Its classical radial limit is
\begin{equation}\label{mod:classicalradial}
  J_C^1u=u''+6r^{-1}u'+6r^{-2}u.
\end{equation}
We first determine the radial solutions of the linearized equation.
Their exponents depend on $\sigma$. We will use the corresponding
powers in the approximation to the leaf, and solve for a smaller
remainder in the fixed weight $\beta$.
Set $b_*=(6-\sigma)/2$, $k_\sigma=(8+\sigma)/2$, and
\begin{equation}\label{mod:logkernel}
  K_\sigma(s)=2\pi^2c_{7,\sigma}
  \int_{-1}^1\int_{-1}^1
  \frac{\sqrt{1-u^2}\sqrt{1-v^2}}
  {(2\cosh s-u-v)^{k_\sigma}}\,\dd u\,\dd v.
\end{equation}
This kernel is even and positive. If $X_1\in C$ has radius one, put
\begin{equation*}
  V_\sigma=c_{7,\sigma}\int_C
  \frac{1-\nu_C(X_1)\cdot\nu_C(Y)}{|X_1-Y|^{8+\sigma}}\dd S_Y>0.
\end{equation*}
The change $r=e^\ell$, $U(\ell)=u(e^\ell)$ gives
\begin{equation*}
  r^\tau\widehat\J_C^\sigma u(r)
  =\PV\int_\R K_\sigma(s)e^{b_*s}
  [U(\log r+s)-U(\log r)]\dd s+V_\sigma U(\log r).
\end{equation*}
This is \cite[Section 11.3, (11.7)]{DDPW}, specialized to dimensions $(4,4)$ and opening one, after the change $\rho=e^s$.
The radius here is the full Euclidean radius in $\R^8$.
The principal value is equivalently computed by adding the integrands at $s$ and $-s$:
their sum vanishes quadratically, and the difference between
geometric and logarithmic truncations tends to zero. This also
justifies substitution of complex powers below.
For $\sigma$ close to one,
\begin{equation}\label{mod:kernelbounds}
  K_\sigma(s)\asymp\eps|s|^{-2-\sigma}\quad(0<|s|<1),
  \qquad K_\sigma(s)\le C\eps e^{-k_\sigma|s|}\quad(|s|\ge1).
\end{equation}
The first estimate follows by comparing the angular integral near its
singularity with $\int_{\R^6}(s^2+|z|^2)^{-k_\sigma}\dd z$; the second
follows directly from \eqref{mod:logkernel}. More precisely, near
$u=v=1$ the angular metric is uniformly equivalent to the Euclidean
metric and $2\cosh s-u-v$ is comparable to $s^2+|z|^2$. Rescaling
$z=|s|\widetilde z$ gives $|s|^{6-2k_\sigma}=|s|^{-2-\sigma}$.
The complementary angular region is bounded. Derivatives of the
kernel satisfy the corresponding bounds away from zero.
Tangent-plane normalization yields
\begin{equation*}
  \int_{|s|<1}s^2K_\sigma(s)\dd s\longrightarrow2,
  \qquad \int_{|s|<1}s^4K_\sigma(s)\dd s\longrightarrow0,
  \qquad V_\sigma\longrightarrow6.
\end{equation*}
These are \cite[Appendix A, Lemmas A.2--A.4]{DDPW}, applied on a
unit conical annulus, with normalization
$c_{7,\sigma}|S^6|\int_0^1s^{-\sigma}\dd s/14=1$.
The complementary integrals are $O(\eps)$ by
\eqref{mod:kernelbounds}; hence the local limits apply despite the
vertex and the unbounded cone.

\begin{definition}\label{mod:symboldefinition}
  The radial indicial symbol is
  \begin{equation}\label{mod:symbol}
    I_\sigma(z)=V_\sigma+\int_\R K_\sigma(s)
    \big[\cosh((b_*-z)s)-\cosh(b_*s)\big]\dd s,
    \qquad-\tau<\Re z<7.
  \end{equation}
  Thus $\widehat\J_C^\sigma r^{-z}=I_\sigma(z)r^{-z-\tau}$.
\end{definition}

For real $0<z<6-\sigma$, this is exactly the radial-power formula
of \cite[Section 11.3]{DDPW}:
\begin{equation*}
  I_\sigma(z)=c_{7,\sigma}
  \big[A_0(4,4,\sigma)^2-C(4,4,\sigma,z)\big],
\end{equation*}
where $A_0$ and $C$ have the meanings assigned there. In particular
$V_\sigma=c_{7,\sigma}A_0(4,4,\sigma)^2$. We retain this real
radial computation; the larger complex strip, uniform estimates and
Mellin remainder needed below are established here.

\begin{lemma}\label{mod:roots}
  The symbols converge, with their complex derivatives, locally uniformly
  on $-2<\Re z<7$ to $I_1(z)=(z-2)(z-3)$. For $\sigma$ sufficiently close
  to one their only real zeros are two simple roots $m,n$, with
  \begin{equation}\label{mod:rootrelations}
    m\longrightarrow2,\quad n\longrightarrow3,\quad m+n=6-\sigma,
    \qquad I_\sigma(b)<0\quad(m<b<n).
  \end{equation}
  On each fixed closed substrip of $-2<\Re z<7$,
  \begin{equation}\label{mod:highfrequency}
    \Re I_\sigma(y+i\xi)\le C-c|\xi|^\tau\qquad(|\xi|\ge1).
  \end{equation}
\end{lemma}

\begin{proof}
  The integrand in \eqref{mod:symbol} vanishes quadratically at zero;
  its tails are integrable precisely when
  $|b_*-\Re z|<k_\sigma$, which is the stated strip. Every compact
  subset of $-2<\Re z<7$ lies in this strip once $\sigma$ is sufficiently
  close to one. Subtracting
  the tangent affine function in the geometric integral, the quadratic
  term converges to \eqref{mod:classicalradial}, because
  $c_{7,\sigma}|S^6|\int_0^1s^{-\sigma}\dd s/14=1$. The remainder and
  the remote integral tend to zero uniformly on compact sets: the latter
  has an exponentially integrable majorant times $\eps$, and the local
  Taylor remainder has higher order than the quadratic term. Cauchy's
  formula on slightly larger compact sets gives convergence of complex
  derivatives. The symmetry $I_\sigma(z)=I_\sigma(2b_*-z)$ is
  \cite[Remark 11.1]{DDPW}. Strict convexity on the real interval,
  $I_\sigma''(y)=\int s^2K_\sigma(s)\cosh((b_*-y)s)\dd s>0$, proves
  \eqref{mod:rootrelations}. Indeed, the argument principle gives one
  simple zero near each classical root. Conjugation symmetry makes both
  zeros real, strict convexity excludes additional real zeros, and
  reflection about $b_*$ gives their sum.
  Finally,
  \begin{equation*}
    \Re I_\sigma(y+i\xi)=I_\sigma(y)-\int_\R K_\sigma(s)
    \cosh((b_*-y)s)(1-\cos\xi s)\dd s.
  \end{equation*}
  Restricting to $|s|<|\xi|^{-1}$ gives the lower bound
  $c\eps\xi^2\int_0^{|\xi|^{-1}}s^{-\sigma}\dd s=c|\xi|^\tau$.
\end{proof}

In particular there are no complex zeros on a line $m<\Re z<n$.
For $m<b<n$ and $u=r^{-b}v(\log r)$ we obtain
\begin{equation*}
  r^{b+\tau}\widehat\J_C^\sigma u
  =\mathcal A_{\sigma,b}v-\kappa_b v,
  \quad \kappa_b=-I_\sigma(b)>0,
  \quad \mathcal A_{\sigma,b}v
  =\PV\int_\R K_\sigma(s)e^{(b_*-b)s}[v(\ell+s)-v(\ell)]\dd s.
\end{equation*}
Positivity of this kernel gives, for zero exterior Dirichlet data on
any interval,
\begin{equation*}
  \|v\|_\infty\le\kappa_b^{-1}
  \| (\mathcal A_{\sigma,b}-\kappa_b)v\|_\infty.
\end{equation*}
At a positive maximum every increment in $\mathcal A_{\sigma,b}$
is nonpositive, and comparison with the constant
$\kappa_b^{-1}\|g\|_\infty$ proves the upper estimate for a source
$g$. Apply the same argument to $-v$ for the lower estimate. The
constant is independent of the interval. This is the comparison
mechanism that will survive after modification across the smooth core.

\begin{lemma}\label{mod:radialestimates}
  For a function cut off near the vertex, of scaled order $\tau+\alpha$
  and weight $r^{-b}$, $-\tau<b<7$,
  \begin{equation*}
    \|\widehat\J_C^\sigma g\|_{C^\alpha_{\rm sc}}
    \le Cr^{-b-\tau}.
  \end{equation*}
  A compact change contributes $O(\eps r^{-8-\sigma})$. If $M_h$ has a
  fixed smooth bounded part and end height $h=O_{C^k_{\rm sc}}(r^{-a_h})$, with
  $a_h$ close to two, then for $2<b<3$, $b+a_h<6$,
  \begin{equation}\label{mod:geometriccomparison}
    \|\widehat\J_{M_h}^\sigma g-\widehat\J_C^\sigma g\|_{C^\alpha_{\rm sc}}
    \le Cr^{-b-a_h-\tau-1}+C\eps r^{-8-\sigma}.
  \end{equation}
\end{lemma}

\begin{proof}
  Split into $\rho<r/2$, $r/2<\rho<2r$, and $\rho>2r$.
  The first and third source integrals are bounded by
  $C\eps r^{-8-\sigma}\int_1^{r/2}\rho^{6-b}\dd\rho$ and
  $C\eps\int_{2r}^\infty\rho^{-2-\sigma-b}\dd\rho$.
  In the middle region subtract the tangent affine increment and rescale.
  The normalized integral of the remaining second difference is bounded
  uniformly by the scaled $C^{\tau+\alpha}$ norm. These three regions
  show why the admissible weights have endpoints $-\tau$ and seven;
  neither restriction can be dropped from the remote estimates.
  For \eqref{mod:geometriccomparison} the local relative change is
  $O(r^{-a_h-1})$; the inner geometric error is bounded by
  $C\eps r^{-8-\sigma}\int_1^{r/2}\rho^{5-b-a_h}\dd\rho$.
  Keep the quadratic cancellation of $1-\nu(X)\cdot\nu(Y)$ in estimating
  the potential, rather than estimating the two normals separately.
  The displacement of the evaluation point changes the remote kernel by
  its first derivative times $O(r^{-a_h})$, and gives the same relative
  factor $r^{-a_h-1}$. The outer radial moments converge with this factor.
  Changes in the fixed completion have bounded mass and distance
  comparable to $r$, explaining their faster tail. These integrals give
  the asserted powers. Differentiating
  remote kernels and applying the local estimate with exponent $\tau+\alpha$ gives the
  scaled H\"older bounds.
\end{proof}

\subsection{The fractional foliation leaf}\label{mod:leafsection}

Let $M_0\subset\{|v|>|u|\}$ be the invariant classical smoothing of
$C$, normalized by
\begin{equation}\label{mod:classicalleaf}
  h_0(r)=a_0r^{-2}+b_0r^{-3}+O_{C^k_{\rm sc}}(r^{-8}),
  \qquad a_0>0,\quad b_0<0.
\end{equation}
We use the classical smoothing and foliation of \cite{BDG,HS};
the expansion and the sign of $b_0$ are
\cite[Proposition 3.3]{Szekelyhidi}. In particular its dilation field satisfies
\begin{equation}\label{mod:dilationclassical}
  \zeta=X\cdot\nu_{M_0}>0,\qquad J_{M_0}^1\zeta=0,
  \qquad\zeta=3a_0r^{-2}+4b_0r^{-3}+O_{C^k_{\rm sc}}(r^{-8}).
\end{equation}
Fix a large $R_0$. Join the classical core to the end
\begin{equation*}
  h_\sigma^0=a_0r^{-m}+q_0(r),\qquad
  q_0=h_0-a_0r^{-2}=O_{C^k_{\rm sc}}(r^{-3}),\quad r>2R_0,
\end{equation*}
using a fixed interpolation on $R_0<r<2R_0$, and call the resulting
surface $M_\sigma^0$. It converges to $M_0$ on compact sets. Use the
spaces $X_{\sigma,\beta}$, $Y_{\sigma,\beta}$ of the preliminaries,
with seven-dimensional charts and a positive lower bound for $r$.
All scaled derivatives of its end are uniformly small after increasing
$R_0$.

\begin{lemma}\label{mod:leafresidual}
  The normalized curvature $E_\sigma^0=\HH_\sigma(M_\sigma^0)$ satisfies
  \begin{equation}\label{mod:leafresidualbound}
    \|E_\sigma^0\|_{C^\alpha_{\rm sc}}
    \le C(r^{-4-\sigma}+r^{-2m-2-\sigma}+\eps r^{-8-\sigma}),
  \end{equation}
  and
  \begin{equation*}
    \|E_\sigma^0\|_{Y_{\sigma,\beta}}\longrightarrow0.
  \end{equation*}
\end{lemma}

\begin{proof}
  The leading power is killed by $I_\sigma(m)=0$; its cutoff contributes
  only an integral over a fixed bounded region. Lemma~\ref{mod:radialestimates} applies to $q_0$.
  Local quadratic terms have size $r^{-\sigma}(r^{-m-1})^2$; the inner
  remote terms are bounded by
  \[
  C\eps r^{-8-\sigma}\int_{R_0}^{r/2}\rho^{5-2m}\dd\rho
  +C\eps r^{-9-\sigma}\int_{R_0}^{r/2}\rho^{6-2m}\dd\rho.
  \]
  The outer terms give the same power and the bounded part of the leaf gives
  $C\eps r^{-8-\sigma}$. This proves \eqref{mod:leafresidualbound},
  including its H\"older version by the local graph calculus and remote
  kernel differentiation. On $r>R$ its weighted bound is $CR^{\beta-3}$.
  On compact sets normalized curvature converges to the vanishing
  classical curvature by \cite[Appendix A, Lemma A.1]{DDPW}; the remote
  tail is uniformly integrable by the preceding bounds. Choose first
  $R$, then $\sigma$.
\end{proof}

We construct the barriers needed below simultaneously for
$b\in\{b_s,\beta,\beta_1\}$. Choose a positive
smooth radial $f_b$, equal to $r^{-b-2}$ on the end. If $s$ is
meridian arclength and $W(s)$ the orbit-area density, then
\[
B_b(s)=\zeta(s)\int_s^\infty
\frac{\int_0^tW(a)\zeta(a)f_b(a)\dd a}{W(t)\zeta(t)^2}\dd t
\]
is smooth at the collapsing orbit and positive. The radial identity
$J_{M_0}^1(\zeta v)=(W\zeta)^{-1}(W\zeta^2v')'$ gives
\begin{equation*}
  J_{M_0}^1B_b=-f_b,\qquad
  B_b=c_b r^{-b}+O_{C^k_{\rm sc}}(r^{-3}),\qquad
  c_b=[(b-2)(3-b)]^{-1}.
\end{equation*}
Define $B_b$ on $M_\sigma^0$ by the normal identification with $M_0$.
The radial estimates give
\begin{equation}\label{mod:barriererror}
  \widehat\J_{M_\sigma^0}^\sigma B_b
  =c_b I_\sigma(b)r^{-b-\tau}
  +O(r^{-4-\sigma}+r^{-b-m-2-\sigma}+\eps r^{-8-\sigma}).
\end{equation}
The sign on the end follows from $I_\sigma(b)<0$ and on the core
from classical convergence. To make these choices uniform, first fix
an end threshold large enough that the errors in
\eqref{mod:barriererror} are smaller than half the leading negative
term. This is possible because $b<3$ and $b+m<6$ for all three weights.
On the remaining fixed compact set, $f_b$ has a positive minimum
after multiplication by the radius weight; convergence of the
normalized operator preserves that strict sign for $\sigma$ close to
one. Hence
\begin{equation}\label{mod:strictleafbarrier}
  \widehat\J_{M_\sigma^0}^\sigma B_b\le-cr^{-b-\tau},
  \qquad cr^{-b}\le B_b\le Cr^{-b}.
\end{equation}

We record an inverse estimate which will be applied first on the leaf
and then on the completed graph. We use the identity in the proof of
\cite[Corollary B.1]{DDPW}; with a strict supersolution, the resulting zeroth-order term is positive.

\begin{proposition}[An inverse estimate]\label{mod:leafinverse}
  Let $S$ be a smooth hypersurface with the uniform scaled charts of
  Section~\ref{sec:preliminaries}, and let $L$ be its normalized Jacobi
  operator. Assume its weighted remote integrals satisfy the bounds
  required in \eqref{glu:schauder}, its potential is bounded by
  $Cr^{-\tau}$ in scaled $C^\alpha$, and there are positive functions
  $B\asymp r^{-\beta}$, $B_s\asymp r^{-b_s}$ such that
  \[
  LB\le-cr^{-\beta-\tau},\qquad LB_s<0.
  \]
  Then $L$ has a unique inverse satisfying
  \begin{equation}\label{mod:leafinverseestimate}
    \|L^{-1}f\|_{X_{\sigma,\beta}(S)}
    \le C\|f\|_{Y_{\sigma,\beta}(S)}.
  \end{equation}
  The constant is uniform for a family with these uniform bounds, and
  invariant sources have invariant solutions.
\end{proposition}

\begin{proof}
  Write $K$ for the positive symmetric normalized kernel. For a compactly
  supported $v=B\psi$, the algebra in the proof of
  \cite[Corollary B.1, following (B.5)]{DDPW} gives
  \begin{equation}\label{glu:picone}
    -\int_S vLv
    =\frac12\iint_{S\times S}B(X)B(Y)(\psi(X)-\psi(Y))^2K(X,Y)
    \dd S_X\dd S_Y+\int_S\frac{-LB}{B}v^2\dd S.
  \end{equation}
  The last term controls $c\int r^{-\tau}v^2$, hence the $L^2$ norm
  on bounded domains and, with the bounded potential, the fractional
  energy. The zero-exterior problem is solved by the variational
  construction described before \cite[Lemma 7.2]{DDPW}; the local
  cancellation and assumed tails justify the same truncation and limit
  procedure.

  At a positive maximum of $u/B$ the jump increments are nonpositive and
  $L(B(u/B))\le(u/B)LB<0$. Comparison with
  $\pm C\|f\|_YB$ gives $\sup r^\beta|u|\le C\|f\|_Y$,
  uniformly in the exhaustion domain; for weak solutions the same
  comparison follows by testing with the positive part and using
  \eqref{glu:picone}. The localized estimate \eqref{glu:schauder}
  then gives \eqref{mod:leafinverseestimate} on interior charts.
  Compactness on these charts and exhaustion give a global solution;
  every fixed chart is eventually interior, so its uniform bound passes
  to the limit. The remote hypotheses justify passage to $L$.
  For uniqueness, $u/B_s\to0$ for a homogeneous solution in weight
  $\beta$, and the same maximum argument gives $u=0$. Uniqueness
  preserves the symmetries.
\end{proof}

For $S=M_\sigma^0$, the strict barriers are
\eqref{mod:strictleafbarrier}. The potential is bounded by
$Cr^{-\tau}$ by quadratic cancellation of its numerator locally;
its remote part and the weighted solution tails are precisely
Lemma~\ref{mod:radialestimates}. Thus the proposition applies with
$L=\widehat\J_{M_\sigma^0}^\sigma$. Its compact classical limit,
used in transferring the barriers, is
\cite[Appendix A, Lemmas A.2 and A.4]{DDPW}: the difference term
converges to $\Delta_{M_0}$ and the normal potential to $|A_{M_0}|^2$
under our normalization. Their tail assumption is supplied here by the
radial moments, rather than by compact convergence alone.

The normal-graph expansion of normalized curvature is
$E_\sigma^0+Lu+Q_\sigma(u)$. The first variation is
\cite[Proposition B.1, (B.1)]{DDPW}, with the opposite sign and
normalization specified in Section~\ref{sec:preliminaries}; the local
tangent-graph expansion is \cite[Section 8, (8.1) and the proof of Lemma 8.1]{DDPW}.
With the local $C^{\tau+\alpha}$ estimate of the preliminaries and the
radial moments below, these formulas give
\begin{equation}\label{mod:leafnonlinear}
  \|Q_\sigma(u)-Q_\sigma(v)\|_Y
  \le C(\|u\|_X+\|v\|_X)\|u-v\|_X.
\end{equation}
The local Hessian calculation in \cite[Lemmas 8.3--8.5]{DDPW},
localized in our uniformly smooth scaled charts with the $C^{\tau+\alpha}$
estimate of the preliminaries, gives the quadratic scale
$r^{-2\beta-2-\sigma}$. The new inner remote moment is
\[
\eps r^{-8-\sigma}\int_1^{r/2}\rho^{5-2\beta}\dd\rho
\le Cr^{-2\beta-2-\sigma},
\]
using $2\beta+1<7$; the outer moment converges. Remote kernel
center differences give the same scaled H\"older bound. Taylor's
integral formula therefore gives \eqref{mod:leafnonlinear}; multiplying
by $r^{\beta+\tau}$ leaves the bounded factor $r^{-\beta-1}$.

\begin{proposition}\label{mod:leaf}
  For $\sigma<1$ sufficiently close to one there is an invariant
  one-sided fractional minimal hypersurface $M_\sigma$, converging
  smoothly on compact sets to $M_0$, whose end height satisfies
  \begin{equation}\label{mod:leafleading}
    h_\sigma=a_0r^{-m}+O_{C^{\tau+\alpha}_{\rm sc}}(r^{-\beta}).
  \end{equation}
  Its dilation field is positive, and its dilates foliate the chosen side of $C$.
\end{proposition}

\begin{proof}
  Apply contraction to $u=-L^{-1}(E_\sigma^0+Q_\sigma(u))$ in a ball of
  radius $2C\|E_\sigma^0\|_Y=o(1)$. Equations
  \eqref{mod:leafinverseestimate} and \eqref{mod:leafnonlinear} show that
  this map preserves the ball and has Lipschitz constant less than
  one half once $\sigma$ is sufficiently close to one. Its fixed point
  solves the exact fractional equation, rather than only an asymptotic
  equation on the end. Since $\beta>m$, the leading
  coefficient remains $a_0$. Interior bootstrapping gives compact smooth
  convergence. On the end the dilation field equals
  $(1+m)a_0r^{-m}+O(r^{-\beta})>0$; on the remaining compact set it is
  positive by \eqref{mod:dilationclassical}. The same two-region argument
  keeps the leaf on its prescribed side. Its radial projection is the
  classical radial projection perturbed on a compact set and with the
  same conical limit on the end; positivity of the dilation field makes
  it a local diffeomorphism. Properness and degree one then give the
  foliation by dilates.
\end{proof}

We also need these barriers on the exact leaf. Define $B_b$ on
$M_\sigma$ by its normal identification with $M_\sigma^0$, for $b\in\{b_s,\beta,\beta_1\}$. Denote this
last correction by $u_\sigma$ and put
$\delta_\sigma=\|u_\sigma\|_{X_{\sigma,\beta}}\to0$. The local change
of the Jacobi kernel and potential acting on $B_b$ is bounded by
\[
C\delta_\sigma r^{-b-\beta-\tau-1}.
\]
For the inner remote part, the geometric comparison of
Lemma~\ref{mod:radialestimates} gives
\[
C\eps\delta_\sigma r^{-8-\sigma}
\int_1^{r/2}\rho^{5-b-\beta}\dd\rho
+C\eps\delta_\sigma r^{-8-\sigma}.
\]
The integral has the same local power because $b+\beta<6$; the outer
moment converges. These estimates hold in scaled $C^\alpha$. They
are small relative to the strict source, on the end and on the fixed
bounded part, respectively. We obtain a positive barrier $B_{{\rm leaf},\sigma}$ with
\begin{equation}\label{mod:exactleafbarrier}
  \widehat\J_{M_\sigma}^\sigma B_{{\rm leaf},\sigma}
  \le-cr^{-b-\tau},\qquad
  B_{{\rm leaf},\sigma}=c_br^{-b}+O_{C^k_{\rm sc}}(r^{-3}).
\end{equation}
The normal change of variables preserves this expansion, and
the constants are uniform for the three selected weights.

\subsection{The two-term expansion and the leaf parameter}\label{mod:mellinsection}

The interchange $(u,v)\mapsto(v,u)$ exchanges the two sets and sends
a radial height $h$ to $-h$. Therefore the curvature pulled back with its evaluation point
is odd in $h$: its quadratic term vanishes. The bounded part of the leaf contributes a remote integral. Consequently,
\begin{equation}\label{mod:cubicremainder}
  Q_C^\sigma(h_\sigma)=O_{C^k_{\rm sc}}(r^{-3m-3-\sigma})
  +O_{C^k_{\rm sc}}(\eps r^{-8-\sigma}).
\end{equation}
Here and below $k$ is any fixed finite integer. We give the remote
calculation, since it determines the available Mellin strip. For a
source radius $\rho$ and separation $D$, the joint pullback of source
and evaluation point bounds the third variation by
\[
C\eps D^{-8-\sigma}\rho
\left(\frac{H_\rho}{\rho}+\frac{H_r}{D}\right)^3\dd S_Y,
\qquad H_s=Cs^{-m}.
\]
To obtain this bound, prolong the normal displacement across a tube
whose width is comparable to the source radius. Each source variation
costs $H_\rho/\rho$, while moving the evaluation point costs $H_r/D$.
Integration across the tube supplies its width $\rho$. This is a
calculation in which both the source and the evaluation point move.
It includes all the mixed terms. Choosing this extension to commute with $(u,v)\mapsto(v,u)$
preserves oddness on each annulus, so the
second variation integrates to zero before the cubic bound is used.
Outside the normal tube only
evaluation-point variations remain; subtraction of the tangent
half-space bounds them by $Cr^{-\sigma}(H_r/r)^3$.
On comparable local radii the same subtraction gives
$C\eps r^{-1}(H_r/r)^3\int_0^{cr}D^{-\sigma}\dd D$, with the same
uniform bound. For
$\rho<r/2$, the cubic pullback is bounded by
\[
C\eps r^{-8-\sigma}\int_{R_0}^{r/2}\rho^{4-3m}\dd\rho
+C\eps\sum_{j=1}^3r^{-jm-8-\sigma-j}
\int_{R_0}^{r/2}\rho^{4+j-(3-j)m}\dd\rho.
\]
The first radial integral is bounded because $3m>5$; the other
three have size $O(r^{-3m-3-\sigma})$. For $\rho>2r$ the four terms are
\[
C\eps\sum_{j=0}^3r^{-jm}
\int_{2r}^\infty\rho^{-(3-j)m-4-\sigma}\dd\rho,
\]
and give the local power. On comparable annuli, tangent affine
subtraction gives $Cr^{-\sigma}(h/r)^3$.
The fixed bounded part of the leaf is not required to satisfy the interchange cancellation
in conical coordinates. The change in the corresponding set has bounded ambient
volume, hence contributes $O(\eps r^{-8-\sigma})$, including
derivatives of the remote evaluation kernel. The contribution $O(\eps r^{-8-\sigma})$ from this bounded part
determines the Mellin strip below.

\begin{lemma}\label{mod:annularbootstrap}
  For $H_R(\rho)=R^{-1}h_\sigma(R\rho)$,
  \begin{equation*}
    \|H_R\|_{C^k(3/4<\rho<3/2)}\le C_kR^{-m-1}.
  \end{equation*}
  In particular \eqref{mod:cubicremainder} holds after finitely many
  applications of $r\partial_r$.
\end{lemma}

\begin{proof}
  The initial bound is \eqref{mod:leafleading} after rescaling. Difference
  quotients of the graph equation satisfy uniformly elliptic equations
  with the kernel of order $\tau$ used in Proposition~\ref{mod:leafinverse}.
  Their coefficient differences have the current graph bound. Their
  remote terms have the same bound by radial integration, since
  $m+1<7$; the fixed bounded part contributes $O(\eps R^{-7})$.
  Here one applies the equation to the difference between the leaf and
  the conical graph on the normalized annulus. Freezing the tangent
  density gives an even principal kernel; terms with an additional
  increment are absorbed by the lower-order interpolation estimate. A
  spatial difference quotient therefore has a right-hand side controlled
  by the already bounded graph difference and complementary integrals.
  The local Schauder estimate on a smaller chart gains a derivative
  without changing $R^{-m-1}$. Repeat on a nested finite collection of
  annuli. At each step the coefficient derivatives involve only derivatives of $H_R$
  controlled at earlier steps, and differentiating a remote kernel at
  the evaluation point costs $R^{-1}$ before rescaling. A finite covering
  then gives the constants $C_k$ independently of $R$ and $\sigma$.
  Iteration proves the assertion.
  Differentiating the preceding cubic integrals preserves their powers.
\end{proof}

Set $v=\chi(r)(h_\sigma-a_0r^{-m})$, with a fixed cutoff vanishing
near zero, and $f=\widehat\J_C^\sigma v$. Exact minimality and
\eqref{mod:cubicremainder} imply, for every fixed $\kappa<7$,
\begin{equation}\label{mod:mellinsource}
  |(r\partial_r)^jf(r)|\le C_{j,\kappa}
  \begin{cases}1,&r\le1,\\r^{-\tau-\kappa},&r\ge1.\end{cases}
\end{equation}
Indeed $3m+2>7$, while the integral over the bounded part decays as
$r^{-\tau-7}$.
The cutoff and the modification of the bounded part contribute the same
remote power. On the intermediate fixed annulus all derivatives are
uniformly bounded. Near zero $v$ vanishes, so the source is an integral
over points at a fixed positive distance and is smooth with bounded
logarithmic derivatives. These observations prove both ends of
\eqref{mod:mellinsource}; they will permit a contour shift with
derivatives of its remainder.

\begin{proposition}\label{mod:mellin}
  There are continuous coefficients $b_\sigma\to b_0<0$ such that
  \begin{equation}\label{mod:leaftwopowers}
    h_\sigma=a_0r^{-m}+b_\sigma r^{-n}+N_\sigma(r),
    \qquad |(r\partial_r)^jN_\sigma|\le C_jr^{-6}.
  \end{equation}
\end{proposition}

\begin{proof}
  Define
  \[
  \mathcal F_\sigma(z)=\int_0^\infty f(r)r^{z+\tau}\frac{\dd r}{r},
  \qquad\mathcal V_\sigma(z)=\int_0^\infty v(r)r^z\frac{\dd r}{r}.
  \]
  The first is holomorphic on $-\tau<\Re z<7$ by
  \eqref{mod:mellinsource}. The second exists for $\Re z<\beta$, since
  $v$ vanishes near zero and has weight $\beta$ on the end. Choose a
  fixed $c$ with $m<c<\beta$, separated from both endpoints. On this line,
  \begin{equation}\label{mod:mellinidentity}
    \mathcal F_\sigma(z)=I_\sigma(z)\mathcal V_\sigma(z).
  \end{equation}
  To justify Fubini, first truncate the logarithmic integral symmetrically.
  Add the logarithmic integrands at $s$ and $-s$ near zero; their weighted integral is bounded
  by $Cs^2$ times the integrable weighted second-derivative norm of $v$.
  For $|s|>1$, the exponential kernel and the weights on the line $c$
  give an integrable bound. The change of variable $r\mapsto re^s$ is
  then legitimate, and removal of the truncations gives
  \eqref{mod:mellinidentity}. The only zero in
  $c\le\Re z\le6$ is $n$: \eqref{mod:highfrequency} excludes large
  imaginary parts, and locally uniform convergence to $(z-2)(z-3)$
  excludes any additional zeros in the remaining compact rectangle.
  This also gives a uniform positive lower bound for $|I_\sigma|$ on
  the lines $\Re z=c$ and $\Re z=6$ at bounded frequencies. At large
  frequencies \eqref{mod:highfrequency} gives
  $|I_\sigma(y+i\xi)|\ge c|\xi|^\tau$ throughout this strip.
  Moreover, \eqref{mod:mellinsource} gives
  \begin{equation}\label{mod:mellinFourierdecay}
    \sup_{c\le y\le6}|\mathcal F_\sigma(y+i\xi)|
    \le C_j(1+|\xi|)^{-j}.
  \end{equation}
  To see this, integrate by parts $j$ times in the Fourier transform of
  $e^{(y+\tau)\ell}f(e^\ell)$; its derivatives are uniformly integrable,
  using any $6<\kappa<7$ in \eqref{mod:mellinsource}.
  In fact their majorants decay like $e^{-(\kappa-y)\ell}$ as
  $\ell\to+\infty$ and $e^{(y+\tau)\ell}$ as $\ell\to-\infty$,
  uniformly for $c\le y\le6$. Fourier inversion for
  $e^{c\ell}v(e^\ell)$ now gives Mellin inversion on the first line.
  Shift through a rectangle to the line six. The horizontal integrals
  tend to zero by \eqref{mod:mellinFourierdecay} and the symbol bound.
  The simple pole at $n$ is the only one crossed. The resulting formula is
  \[
  v(r)=-\frac{\mathcal F_\sigma(n)}{I_\sigma'(n)}r^{-n}
  +\frac1{2\pi i}\int_{6-i\infty}^{6+i\infty}
  r^{-z}\frac{\mathcal F_\sigma(z)}{I_\sigma(z)}\dd z.
  \]
  Thus
  \begin{equation}\label{mod:residuecoefficient}
    b_\sigma=-\mathcal F_\sigma(n)/I_\sigma'(n).
  \end{equation}
  The minus sign is the residue sign for a rightward shift with inversion
  factor $r^{-z}$. Each application of $r\partial_r$ multiplies the
  integrand on the line six by $-z$. Choosing the integration-by-parts
  order in \eqref{mod:mellinFourierdecay} larger than the prescribed
  derivative order makes this integral absolutely convergent, and gives
  the same $O(r^{-6})$ bound for every fixed logarithmic derivative.
  Compact convergence of the leaves
  and a common integrable majorant from \eqref{mod:mellinsource} imply
  continuity and convergence of \eqref{mod:residuecoefficient}. At
  $\sigma=1$, \eqref{mod:classicalleaf} identifies the limit with $b_0$.
\end{proof}

For $\lambda>0$ use $\lambda M_\sigma$; for $\lambda<0$ use
$|\lambda|S(M_\sigma)$, where $S(u,v)=(v,u)$. Label these leaves by
$A=a_0\operatorname{sgn}(\lambda)|\lambda|^p$, and set
\begin{equation*}
  p=1+m,\quad d=n-m,\quad q=\frac{1+n}{1+m},\quad
  \phi_q(A)=\operatorname{sgn}(A)|A|^q,\quad b=b_\sigma a_0^{-q}.
\end{equation*}
Here $p\to3$, $d\to1$, $q\to4/3$. Their conical representation is
\begin{equation*}
  \Psi(r,A)=Ar^{-m}+b\phi_q(A)r^{-n}+N(r,A),
  \qquad r\ge K|A|^{1/p}.
\end{equation*}
The complete leaf is always the dilated smooth hypersurface, including
the region near its collapsed orbit. With $\mu_N=7/p$ its remainder satisfies
\begin{equation*}
  |\partial_A^j(r\partial_r)^kN(r,A)|
  \le C_{j,k}|A|^{\mu_N-j}r^{-6}\qquad(A\ne0).
\end{equation*}
Indeed each label derivative lowers the dilation degree by $p$, and
\eqref{mod:leaftwopowers} gives degree seven. Explicitly, the scaled
remainder is $\operatorname{sgn}(A)|\lambda|N_\sigma(r/|\lambda|)$;
differentiation with respect to $A$ produces linear combinations of
$(r\partial_r)^jN_\sigma$ and factors $|\lambda|^{-p}$. Since $\mu_N\to7/3$,
$N$ extends through zero with two vanishing label derivatives, and its
second derivative is uniformly $C^\eta$, $\eta=1/8$, on bounded label
intervals. For comparable labels at separation $\delta$ smaller than
their size, integration of the third-derivative bound gives
$Cr^{-6}\delta |A|^{\mu_N-3}\le C r^{-6}\delta^\eta$.
If the separation is comparable to the larger size, use
$|\partial_A^2N|\le Cr^{-6}|A|^{\mu_N-2}$ and its vanishing at zero.
The positive gap $\mu_N-2-\eta$ makes both estimates uniform near
$\sigma=1$. The same argument applies after any fixed radial
derivatives. Third label derivatives are used only away from zero.
The explicit $\phi_q$ term is kept separate and will be replaced by a
smooth cylindrical field in the inner region. This explains the contour
line six: a remainder of order $r^{-L}$ has label degree $(1+L)/p$,
so two label derivatives require $L>2p-1$, whereas compact mass limits
the Mellin strip to $L<7$.

\section{Jacobi fields on the cylindrical cone}\label{mod:cylindersection}

Write $\Sigma=C\times\R_t$ and $d=n-m$. We use the normalization
\begin{equation*}
  \widehat\J_\Sigma^\sigma=\kappa_\sigma\J_\Sigma^\sigma
\end{equation*}
fixed in \eqref{eq:normalization}.

We next solve the linearized equation on $C\times\R_t$ with prescribed
behavior for $|t|/r$ large. Kelvin inversion of $r^{-m}$ gives positive
solutions, and integration in the vertical variable gives the functions
$H_{\sigma,q}$ below. We then compute their first two asymptotic terms
by the ordinary differential equation they satisfy.

There is also an exact identity for functions which are affine in $t$:
\begin{equation*}
  \widehat\J_\Sigma^\sigma[t\,a(r)]
  =t\widehat\J_C^\sigma a(r).
\end{equation*}
For $a=r^{-b}$ it holds whenever $-\sigma<b<7$. To see this,
write $t'=t+z$; the additional term is odd in $z$ and integrates to
zero. Its absolute outer radial tail is bounded by
$C\int^\infty\rho^{-1-\sigma-b}\dd\rho$ and its inner tail by
$C\int_0\rho^{6-b}\dd\rho$, so the two stated restrictions justify
the integration. In particular, $tr^{-m}$ is an exact cylindrical Jacobi
field; it is the inner mode associated with the principal leaf label.

The exact power is essential: expanding $r^{-m_\sigma}$ about $r^{-2}$
produces logarithms which are not uniformly small on the end.

\subsection{Kelvin covariance}

\begin{proposition}\label{mod:kelvin}
  Let $\mathcal Ku(X)=\abs X^{\tau-8}u(X/\abs X^2)$. Whenever the
  principal values and tails involved are defined,
  \begin{equation}\label{mod:kelvinidentity}
    \J_\Sigma^\sigma(\mathcal Ku)(X)
    =\abs X^{-8-\tau}
    (\J_\Sigma^\sigma u)(X/\abs X^2).
  \end{equation}
\end{proposition}

\begin{proof}
  Inversion preserves the cone $\Sigma$, which has dimension $N=8$,
  changes its surface measure by the factor $\abs Y^{-2N}$, and satisfies
  \[
  \left|\frac X{\abs X^2}-\frac Y{\abs Y^2}\right|
  =\frac{\abs{X-Y}}{\abs X\abs Y}.
  \]
  When changing variables in the difference operator, a possible extra
  term is a multiple of the operator applied to
  $w(X)=\abs X^{-(N-\tau)}$, and we show that this term vanishes. At
  $\abs X=1$, inversion in $Y$ maps
  $(\abs Y^{-(N-\tau)}-1)\abs{X-Y}^{-N-\tau}\dd S_Y$ to its negative. The
  integral converges at the vertex and at infinity, and near $X$ we use a
  symmetric principal value. The image of the small excised ball differs
  from it by an annulus of thickness $O(\delta^2)$, whose uncancelled
  first-order contribution is $O(\delta^{2-\tau})$ and tends to zero.
  Hence the principal value vanishes. This argument applies at every point
  of unit length, and homogeneity extends it to all $X$.

  The change of variables therefore gives the covariance of the difference
  operator without an extra term. The normal of a cone is unchanged under
  inversion, since it is perpendicular to the radial direction, and the
  Jacobi potential is homogeneous of degree $-\tau$, so that
  $V(X/\abs X^2)=\abs X^{2\tau}V(X)$. This is exactly the remaining factor
  in \eqref{mod:kelvinidentity}.
\end{proof}

Since $r^{-m}$ is a Jacobi field and $(8-\tau)-2m=1+n-m=1+d$, Kelvin
inversion and vertical translations produce the positive Jacobi fields
\begin{equation}\label{mod:positivekernels}
  \frac{r^{-m}}{(r^2+(t-s)^2)^{(1+d)/2}},\qquad s\in\R.
\end{equation}
We normalize the corresponding one-dimensional kernel by
\begin{equation*}
  P_d(r,z)=c_d\frac{r^d}{(r^2+z^2)^{(1+d)/2}},\qquad
  c_d=\frac{\Gamma((1+d)/2)}{\sqrt\pi\Gamma(d/2)},
\end{equation*}
so that, by the beta integral, $\int_\R P_d(r,z)\dd z=1$.

\subsection{Odd boundary powers}

\begin{proposition}\label{mod:poisson}
  For $-1<q<1+d$, define
  \begin{equation*}
    H_{\sigma,q}(r,t)=r^{-n}\lim_{L\to\infty}
    \int_{-L}^L P_d(r,t-s)\phi_q(s)\dd s,\qquad
    \phi_q(s)=\mathop{\rm sgn}(s)\abs s^q.
  \end{equation*}
  The symmetric limit exists, and $H_{\sigma,q}$ is smooth for $r>0$, odd
  in $t$ and homogeneous of degree $q-n$. It satisfies
  \begin{equation}\label{mod:poissonharmonic}
    \widehat\J_\Sigma^\sigma H_{\sigma,q}=0,
  \end{equation}
  and its boundary coefficient equals one: for $t\ne0$,
  $r^nH_{\sigma,q}(r,t)\to\phi_q(t)$ as $r\downarrow0$.
\end{proposition}

\begin{proof}
  Add the integrands at $s$ and $-s$. On compact subsets of $\{r>0\}$,
  their sum behaves at infinity like
  \[
  [P_d(r,t-s)-P_d(r,t+s)]s^q=O(s^{q-d-2}),
  \]
  which is integrable when $q<1+d$, while the condition $q>-1$ gives local
  integrability at zero. Derivatives on compact subsets have the same or
  better decay, so the limit is smooth. Oddness and homogeneity follow by
  changes of variables, and the approximate identity property, together
  with the control of the sum of the two tails, gives the boundary coefficient
  away from $t=0$.

  To prove \eqref{mod:poissonharmonic}, we first keep $L$ finite. The
  positive kernels \eqref{mod:positivekernels} are Jacobi fields, and so is
  their truncated superposition, by Fubini's theorem: we use symmetric
  second differences near the evaluation point and absolute integrability
  elsewhere. We must then pass to the limit $L\to\infty$ in the tail of
  the Jacobi operator, not only on compact sets. This representation as an integral over $s>0$
  gives the bound, uniform in $L$,
  \begin{equation*}
    \abs{H_{\sigma,q}^{\,L}(r,t)}
    \le Cr^{-n}(r+\abs t)^q.
  \end{equation*}
  For $\abs t\le Cr$, this follows by rescaling by $r$ and estimating the displayed difference of kernels. For $\abs t>Cr$, we split the integral at $\abs t/2$
  and $2\abs t$. On the middle interval we use that the kernel has mass
  one and that $s^q\asymp\abs t^q$, while on the other two intervals we
  use the difference of the kernels and obtain the bound
  $Cr^d\abs t^{q-d}\le C\abs t^q$. The argument allows empty intervals and
  is uniform for the truncated integrals.

  Near the axis $r=0$ the majorant is integrable against $r^6\dd r$ since
  $n<7$; at the origin its degree of homogeneity is $q-n>-8$; and at
  infinity it belongs to the tail space of the Jacobi operator because
  $q-n<\tau$. These inequalities hold throughout the parameter range of
  the proposition for $\sigma$ near one. More explicitly, in polar
  coordinates on the eight-dimensional cylinder, the integral near its
  vertex has radial factor $R^{7+q-n}\dd R$, while the weighted integral
  at infinity has factor $R^{q-n-\tau-1}\dd R$. The angular integral at
  the vertical axis is finite because $n<7$. These estimates provide
  one majorant for the whole sequence of truncated fields, rather than
  only for their restriction to a compact chart.
  Dominated convergence therefore
  applies outside a small fixed ball about the evaluation point, and
  inside that ball we use smooth convergence together with the
  cancellation of second differences. This proves
  \eqref{mod:poissonharmonic}.
\end{proof}

\begin{corollary}\label{mod:centralmodebounds}
  For fixed $M<\infty$ and all fixed finite $i,j\ge0$,
  \[
  \abs{\partial_r^i\partial_t^jH_{\sigma,q}(r,t)}
  \le C_{i,j,M}r^{q-n-i-j},
  \qquad \abs t\le Mr.
  \]
  In particular, $\abs{H_{\sigma,q}(r,t)}\le C_M\abs t\,r^{q-n-1}$ in this
  region. The constants are uniform when $q$ ranges in either matching
  neighborhood and $\sigma$ is close to one.
\end{corollary}

\begin{proof}
  At $r=1$, the integral over $s>0$ of the difference of kernels and its finite derivatives converge
  uniformly on $\abs t\le M$, by the bounds in
  Proposition~\ref{mod:poisson}. Scaling $(r,t)$ by $r$ and using
  homogeneity gives the first estimate. Oddness gives $H_{\sigma,q}(r,0)=0$,
  and integrating the first vertical derivative proves the last estimate.
\end{proof}

\subsection{The second asymptotic term}

Set $F=r^nH_{\sigma,q}$. Direct differentiation of the Poisson kernel
gives
\begin{equation}\label{mod:extensionequation}
  F_{rr}+\frac{1-d}{r}F_r+F_{tt}=0,
\end{equation}
a local differential equation satisfied by this particular family of
cylindrical Jacobi fields. For $t>0$ we write $F(r,t)=t^q f(r/t)$, which
leads to
\begin{equation}\label{mod:frobenius}
  (1+x^2)f''+\left(\frac{1-d}{x}-2(q-1)x\right)f'
  +q(q-1)f=0.
\end{equation}
The indicial exponents of this equation at $x=0$ are $0$ and $d$. Since
$0<d<2$, the corresponding solutions have convergent expansions
$f_0(x)=1+O(x^2)$ and $x^df_1(x)=x^d(1+O(x^2))$.
No identification of the full fractional Jacobi operator with the local
operator in \eqref{mod:extensionequation} is being made: exact Jacobi
harmonicity has already been proved by Kelvin inversion and the tail
limit. The local equation is used only for the asymptotic expansion.

\begin{proposition}\label{mod:twotraces}
  Let $q$ range over a fixed small neighborhood of one of the two matching
  exponents $4/3$ and $2/3$. Then, for $\abs t\ge2r$,
  \begin{equation*}
    H_{\sigma,q}(r,t)=r^{-n}\phi_q(t)
    +k(q,d)r^{-m}\phi_{q-d}(t)+\mathcal R_q(r,t),
  \end{equation*}
  where the connection coefficient is defined by
  \begin{equation}\label{mod:connectioncoefficient}k(q,d)=c_d\,\PV\int_\R
    \frac{\phi_q(s)-1}{|1-s|^{1+d}}\dd s.
  \end{equation}
  The singularity at one is taken symmetrically and the integrands on the positive and negative tails are added before taking the limit. Moreover, for all
  fixed finite $i,j\ge0$,
  \begin{equation}\label{mod:twotraceremainder}
    \abs{\partial_r^i\partial_t^j\mathcal R_q(r,t)}
    \le C_{i,j}r^{2-n-i}\abs t^{q-2-j}.
  \end{equation}
\end{proposition}

\begin{proof}
  The coefficient of $f_0$ is one by Proposition~\ref{mod:poisson}.
  At $t=1$, subtract the unit mass of $P_d$ and divide by $r^d$.
  The resulting expression is
  \[
  \frac{F(r,1)-1}{r^d}
  =c_d\lim_{L\to\infty}\int_{-L}^L
  \frac{\phi_q(s)-1}{(r^2+(1-s)^2)^{(1+d)/2}}\dd s.
  \]
  Choose an even cutoff supported near $s=1$. The integral of
  $q(s-1)$ times that cutoff is zero for every $r$, so subtracting it
  does not change this formula. Cancellation of the odd increment at $s=1$ and addition of the
  integrands at $s$ and $-s$ for large $s$ give \eqref{mod:connectioncoefficient} by dominated
  convergence.
  The remaining local numerator is $O((s-1)^2)$, integrable since $d<2$;
  the sum of the tail integrands is
  $O(s^{q-d-2}+s^{-1-d})$, integrable since $q<1+d$ and $d>0$.
  At zero use $q>-1$. These bounds are locally uniform in the parameters
  and also prove continuity of $k$. Parameter differences insert at most
  logarithmic factors in these integrable bounds, so the coefficients
  vary smoothly on either of the compact matching parameter ranges.

  For completeness, the two Frobenius solutions can be obtained directly
  from \eqref{mod:frobenius}. Write
  $f_\rho=x^\rho\sum_{j\ge0}c_{\rho,j}x^{2j}$, with
  $\rho=0$ or $d$ and $c_{\rho,0}=1$. Their coefficients satisfy
  \begin{equation*}
    c_{\rho,j}=-\frac{(\rho+2j-2-q)(\rho+2j-1-q)}
    {(\rho+2j)(\rho+2j-d)}c_{\rho,j-1},\qquad j\ge1.
  \end{equation*}
  The denominators stay uniformly away from zero for the parameters
  under consideration. The coefficient ratios tend to one, so the
  series converge for $|x|<1$, uniformly with all fixed derivatives on
  $|x|\le1/2$. This also covers $d=1$: the two solutions have opposite
  parity and no logarithmic resonance occurs.

  Finally, the convergent Frobenius expansions give
  \[
  F(r,t)=\phi_q(t)\{1+O((r/t)^2)\}
  +k(q,d)r^d\phi_{q-d}(t)\{1+O((r/t)^2)\},
  \]
  and the differentiated series converge uniformly for $\abs{r/t}\le1/2$
  and for parameters in these neighborhoods. After multiplication by
  $r^{-n}$, the first error term gives \eqref{mod:twotraceremainder}, and
  the second is smaller by a factor $(r/\abs t)^d$. Thus the convergent
  local series provide the remainder together with all its differentiated
  bounds. Note that the kernel itself has an infinite second moment, so it
  is the Frobenius expansion that supplies the second-order term.
\end{proof}

For $d=1$ the formulas reduce to
\begin{equation*}
  k(q,1)=-q\cot(\pi q/2),\qquad
  H_{1,q}(r,t)=\frac{\mathop{\rm Im}(r+it)^q}
  {r^3\sin(\pi q/2)},
\end{equation*}
where $\arg(r+it)\in(-\pi/2,\pi/2)$. The explicit harmonic function
has the prescribed boundary coefficient; its Taylor expansion at
$r/|t|=0$ gives the displayed value of $k(q,1)$.

In particular, the two coefficients used below have the limits
\[
k(4/3,1)=\frac4{3\sqrt3},\qquad
k(2/3,1)=-\frac2{3\sqrt3}.
\]
Neither vanishes. Thus the term proportional to $r^{-m}$ cannot be omitted in the
classical limit, and continuity shows that it persists for the exact
fractional exponents. This observation also distinguishes the two
roles of the expansion: $r^{-n}\phi_q(t)$ determines the coefficient of $r^{-n}$ in the leaf
expansion, while $k(q,d)r^{-m}\phi_{q-d}(t)$ changes the coefficient
of $r^{-m}$.
Both are needed even though their relative sizes depend on the region
in which the field is evaluated. All later uses of $k$ require only
its convergent integral, its parameter continuity, and these limits.

\section{The approximate graph and its error}\label{res:section}

We construct an approximate graph from the exact fractional leaves and
the cylindrical Jacobi fields. The separation of tangent-chart and
remote contributions is the one used in \cite[(4.1), Lemma 4.1]{DDPW}.
For the catenoid, its remote two-sheet interaction is evaluated in
\cite[Lemma 4.2, (4.7)]{DDPW}; together with the local curvature term
it determines the ODE \cite[(2.6)]{DDPW}. Here the reference is an
exact fractional leaf, so its whole curvature cancels. The estimates
below control changes of label, the bounded part of each leaf and the
interpolation between the two approximations.
Proposition~\ref{res:assembled} is the resulting counterpart at infinity of
\cite[Proposition 2.1]{DDPW}, in the present weights and geometry.
Throughout this section,
\[
\eps=1-\sigma,\quad \tau=1+\sigma,\quad m=m_\sigma,\quad n=n_\sigma,
\quad d=n-m,\quad p=1+m,
\]
and use the matching exponent $a=5/2$.
We restrict $\sigma$ so that
\begin{equation*}
  2+\alpha/2<\tau+\alpha<2+\alpha<2+\eta.
\end{equation*}
Write $\Sigma=C_{4,4}\times\R_t$, $r=\abs P$ for $(P,t)\in\Sigma$,
and $\mathcal U_r$ for a chart of radius comparable to $r$. Scaled norms
are computed after dilation by $r$; $\mathcal U_r$ need not lie in an
ambient annulus, since the angular region has $\abs t\asymp r^p$.
We use $E_C=\{\abs u>\abs v\}$ and its outward normal $\nu_C$; a normal
height $h$ represents the set $\{s<h\}$. The upward orientation of the
resulting vertical graph is the opposite orientation.

\subsection{The matching surface}\label{res:construction}

Proposition~\ref{mod:mellin} gives
\begin{equation}\label{res:leafexpansion}
  \Psi(r,A)=Ar^{-m}+b\phi_q(A)r^{-n}+N(r,A),\qquad
  q=\frac{1+n}{p},\qquad \phi_j(A)=\operatorname{sgn}(A)\abs A^j,
\end{equation}
where $b=b_\sigma a_0^{-q}$. With $\mu_N=7/p$, its differentiated remainder satisfies
\begin{equation}\label{res:Nderivatives}
  \abs{\partial_A^j(r\partial_r)^kN(r,A)}
  \le C_{j,k}\abs A^{\mu_N-j}r^{-6},\qquad r\ge K\abs A^{1/p}.
\end{equation}
In particular, $\mu_N>2+\eta$, whereas the explicitly retained power
$\phi_q$ is not twice differentiable at zero. The exact cylindrical fields satisfy
\begin{equation*}
  H_j(r,t)=r^{-n}\phi_j(t)+k(j,d)r^{-m}\phi_{j-d}(t)+R_j(r,t),
  \qquad \abs t\ge2r,
\end{equation*}
\begin{equation*}
  \abs{\partial_r^i\partial_t^lR_j}
  \le C_{i,l}r^{2-n-i}\abs t^{j-2-l}.
\end{equation*}
Set
\begin{equation*}
  q_1=2q-d-1,\qquad q_2=3q-2d-2,
  \qquad a_1=bk(q,d),\quad a_2=bqa_1k(q_1,d).
\end{equation*}
Their limits are $q=4/3$, $q_1=2/3$, $q_2=0$. For large $\abs t$ let
\begin{equation*}
  A_{\rm out}(t)=t+a_1\phi_{q-d}(t)+a_2\phi_{q_1-d}(t),
\end{equation*}
and complete this label smoothly and oddly on a fixed bounded interval,
with $A_{\rm out}(t)=t$ near zero. Choose even smooth cutoffs $\chi_A$,
zero near zero and one for large arguments, and $\chi_0$, zero on
$[-1,1]$ and one outside $[-2,2]$. Put
\begin{equation*}
  A_{\rm eff}=t+\chi_A(t/r)(A_{\rm out}(t)-t),
\end{equation*}
\begin{equation}\label{res:innerheight}
  \psi_{\rm in}=\ell+\chi_0(t/r)N(r,A_{\rm eff}),\qquad
  \ell=tr^{-m}+bH_q+bqa_1H_{q_1}.
\end{equation}
Thus
\begin{equation}\label{res:harmonicell}
  \widehat\J_\Sigma^\sigma\ell=0.
\end{equation}
The cutoff of $N$ ensures that the height is smooth
through $t=0$. The term $r^{-n}\phi_q(A)$ in \eqref{res:leafexpansion} cannot be
inserted directly there: its second label derivative is singular.
Instead, it is represented by the smooth cylindrical modes in $\ell$.
The separate completion of $A_{\rm out}$ is needed only on a fixed label
interval; all expansions involving large labels are used outside that
interval.
For $\abs t\ge Cr$ write $S=A_{\rm out}r^{-m}+b\phi_q(A_{\rm out})r^{-n}$
and $W=\ell-S$. Matching the two slow traces gives
\begin{equation}\label{res:Wbounds}
  \abs{\partial_r^i\partial_t^jW}
  \le C_{i,j}\big(r^{2-n-i}\abs t^{q-2-j}+r^{-n-i}\abs t^{q_2-j}\big).
\end{equation}
Indeed, the slow traces of $bH_q$ and $bqa_1H_{q_1}$ have coefficients
$a_1$ and $a_2$, while the first unmatched Taylor power in
$b\phi_q(A_{\rm out})$ is $q_2$.
More explicitly, the first variation of this expression is
$bqa_1\phi_{q_1}(t)$, and its slow trace is precisely
$a_2r^{-m}\phi_{q_1-d}(t)$. This explains why two matching modes, rather
than one, enter \eqref{res:innerheight}. Taylor's formula with the
differentiated completed label gives \eqref{res:Wbounds} for every fixed
number of derivatives. The analytic remainder of $H_q$ gives its first
term, and the first unmatched algebraic power gives its second.

Let $\chi(r,t)=1$ for $\abs t\le r^a$ and $\chi=0$ for
$\abs t\ge2r^a$. Interpolate between $\psi_{\rm in}$ and
$\Psi(r,A_{\rm out}(t))$ using this cutoff. In the angular region retain
the whole exact leaf, including its collapsed orbit, and complete the
fixed compact region by a smooth symmetric entire graph. This defines
$\cP_\sigma$. Since $a<p$, the interface height divided by $r$ is
$O(r^{a-p})$, and its positive $t$-derivative is asymptotic to $r^{-m}$;
the interpolation therefore remains a vertical graph.
In the outer region this graph is described by
\[
A_{\rm out}(t)=F_{\infty,\sigma}(x)
=|x|^pg_\sigma(\vartheta),\qquad
t=A_{\rm out}^{-1}\big(|x|^pg_\sigma(\vartheta)\big),
\]
with the homogeneous function of \eqref{eq:fractionalprofile}.
At $\sigma=1$ this is $F_{\infty,1}=r^3g_1(\vartheta)$,
whose angular function satisfies \eqref{eq:angular-classical} on
the open sector. Thus the power expansion of $A_{\rm out}$ gives the difference between the height of
our approximation and the homogeneous function; the level sets in the
outer region are the exact fractional leaves.
We require no smallness of the curvature on the fixed compact
completion. It defines the set bounded by the approximate graph;
its curvature will be handled by compact matching. The estimates below
are uniform in the end threshold after this completion has been fixed.

\subsection{The bounded part of a leaf}\label{res:endpoint}The movement of the evaluation point must be included in every
variation. For a tubular prolongation $V_h$ of $h\nu_C$, with even
normal cutoff, put $\Phi_\theta=I+\theta V_h$ and
$s_C=\chi_{E_C^c}-\chi_{E_C}$. Then
\begin{equation}\label{res:pullback}
  H^\sigma_{\Phi_\theta(E_C)}(\Phi_\theta(P))
  =\PV\int_{\R^9}s_C(Y)
  \frac{\det(I+\theta DV_h(Y))}
  {\abs{Y-P+\theta(V_h(Y)-V_h(P))}^{9+\sigma}}\,\dd Y.
\end{equation}
We use the tangent-half-space subtraction of \cite[Section 8, before
(8.1)]{DDPW} and the first variation evaluated at the perturbed point
\cite[Proposition B.1, (B.1)]{DDPW}. The latter is stated in arbitrary
ambient dimension. Apply it to compactly supported heights on a smooth
completion away from the vertex, then shrink that completion and remove
the auxiliary cutoffs using the remote bounds below. Using $\chi_{E^c}-\chi_E$ in place of $\chi_E-\chi_{E^c}$ gives
\begin{equation*}
  I'(0)=2\PV\int_\Sigma
  \frac{h(P)\nu_P\cdot\nu_Y-h(Y)}{\abs{Y-P}^{9+\sigma}}\,\dd S_Y
  =-2\J_\Sigma^\sigma h(P).
\end{equation*}
The normalized linearization is $\widehat\J_\Sigma^\sigma$.
Interchanging $u$ and $v$ and applying rotations gives
an isometry fixing $P$, interchanging $E_C\times\R$ with its complement and reversing the normal. For
rotational heights its action sends $V_h$ to $-V_h$. Consequently
\begin{equation}\label{res:F2}
  I(-\theta)=-I(\theta),\qquad I''(0)=0,
\end{equation}
also with equivariant radial and vertical partitions. This cubic
cancellation is special to the present cone. The tangent-graph
calculation of \cite[Lemma 4.1, (4.4)--(4.5)]{DDPW}, differentiated as
in \cite[Section 8]{DDPW}, applies here with eight tangent variables.
After subtracting the first-order Taylor polynomial, the third variation gains $D/r$ and gives
\begin{equation}\label{res:F3}
  \eps\frac{a_h^3}{r}\int_0^{cr}D^{-\sigma}\,\dd D
  \le Cr^{-\sigma}a_h^3,
\end{equation}
where $a_h$ is the scaled height divided by $r$.
Indeed, $V_h(Y)-V_h(P)=DV_h(P)(Y-P)+O(a_hD^2/r)$;
the affine half-space contributes zero for every deformation parameter.
Distances and Jacobians remain uniformly comparable, so the integral
Taylor remainder has the same bound. This verifies the change of
dimension and retains $V_h(P)$ in the kernel. Choose the local cutoff
radially about the moved point. Its displacement derivatives are
supported on $D\asymp r$ and have bounds $Cr^{-j}$; their terms containing $h(P)$ satisfy \eqref{res:farcubic}. The regional partitions preserve
\eqref{res:F2}, and their derivatives cancel on summing the partition.

\begin{lemma}\label{res:localcubic}
  The local normalized nonlinear remainder is bounded by
  $Cr^{-\sigma}a_h^3$, and its difference at $h,\widetilde h$ by
  $Cr^{-\sigma}(a_h+a_{\widetilde h})^2a_{h-\widetilde h}$.
  These bounds hold in scaled $C^\alpha$ with height norms in
  $C^{\tau+\alpha}$ and uniform $C^{2,\eta}$ chart coefficients.
\end{lemma}
\begin{proof}
  Use \eqref{res:F2}--\eqref{res:F3}, the integral Taylor remainder and
  polarization. The normalized difference-quotient estimate of
  Section~\ref{sec:preliminaries} gives the asserted H\"older norm.
\end{proof}

Let $\mathcal Q_{\rm cyl}$ and $\mathcal Q_{\rm rad}$ denote these
normalized remainders, always evaluated at the displaced point; in the
radial expression the vertical variable has been integrated out.
Off the local chart, integration over the tube gives the polarized bound
\begin{equation}\label{res:farcubic}
  C\eps D^{-9-\sigma}\rho
  \left(\frac{H_Y}{\rho}+\frac{H_P}{D}\right)^2
  \left(\frac{W_Y}{\rho}+\frac{W_P}{D}\right)\dd S_Y,
  \qquad D\ge cr.
\end{equation}
Here $H$ bounds height values and $W$ their difference. Indeed,
the tube is $\int_0^{h(Y)}K(X,Y+s\nu_C)J_C(Y,s)\,\dd s$;
its variations differentiate the kernel and $J_C$, not label derivatives
of $h(Y)$. The terms from the complement of this region containing
only heights at $P$ have bound
$Cr^{-\sigma-3}H_P^2W_P$. They cancel when the profile with label
$A_{\rm out}(s)$ and the profile with constant label $A_{\rm out}(t)$
have the same value at the evaluation point.
For example, the pure-$Y$ product is
$C\eps D^{-9-\sigma}H_Y^2W_Y/\rho^2$. The other products contain
one or two factors evaluated at $P$, and are required even though the graph
Jacobian at $Y$ is local. They arise by differentiating the evaluation
point in the kernel. This tube representation permits the use of value
envelopes at labels close to zero, where differentiated two-trace
envelopes would be unavailable.

At the collapsed orbit normal graph coordinates fail. The contribution
that replaces the catenoid's remote two-sheet term is represented by a
signed functional supported in a bounded region, retained before
expanding the graphical end.
For $A\ne0$ set $\lambda=(\abs A/a_0)^{1/p}$ and fix
$\theta=\theta(\rho/\lambda)$, zero on $\rho\le K_0\lambda$ and one
on $\rho\ge2K_0\lambda$. Define
\begin{align*}
  L_A(\varphi)&=-\int_C\theta\Psi\varphi\,\dd S,\\
  Q_A(\varphi)&=-\int_C\theta\int_0^\Psi
  [\varphi(P+s\nu_C)(1-s^2/\rho^2)^3-\varphi(P)]\,\dd s\,\dd S,
  \\
  \mu_A&=(\chi_{E_C}-\chi_{E_A})\,\dd Y-L_A-Q_A.
\end{align*}
Combine the densities before integration. This functional has support
in $B_{C\lambda}$ and satisfies
\begin{equation}\label{res:capscaling}
  \mu_A(\varphi)=\lambda^8\mu_{a_0\operatorname{sgn}A}(\varphi(\lambda\,\cdot)).
\end{equation}
The end contributions cancel exactly because
\[
(L_A+Q_A)(\varphi)
=-\int_C\theta\int_0^{\Psi(\rho,A)}
\varphi(P+s\nu_C)J_C(P,s)\,\dd s\,\dd S_P.
\]
Where $\theta=1$, this is the signed measure of the region between
$C$ and the leaf. After dilation,
the remaining region has bounded volume and bounded tube mass. Thus its
total variation is bounded uniformly for normalized labels. In
particular, \eqref{res:capscaling} is an identity between compact test
functionals with bounded support.
With $E=Y\cdot D_Y$, differentiation gives
\begin{align*}
  \partial_A\mu_A(\varphi)&=\frac{\lambda^8}{pA}
  \mu_{a_0\operatorname{sgn}A}((8+E)\varphi(\lambda\,\cdot)),\\
  \partial_A^2\mu_A(\varphi)&=\frac{\lambda^8}{p^2A^2}
  \mu_{a_0\operatorname{sgn}A}((8+E)(8-p+E)\varphi(\lambda\,\cdot)),
\end{align*}
and hence
\begin{equation}\label{res:captest}
  \abs{\partial_A^j\mu_A(\varphi)}
  \le C\lambda^{8-jp}\sum_{k\le j}\lambda^k
  \norm{D^k\varphi}_{L^\infty(B_{C\lambda})},\qquad j\le2.
\end{equation}
A third differentiation followed by label integration proves H\"older
continuity of the second derivative for exponents below $8/p-2$;
the functional and its first two label derivatives vanish at zero.
These are distributional test bounds.
For completeness, the third Euler polynomial is
$(8+E)(8-p+E)(8-2p+E)$, with prefactor
$\lambda^8/(p^3A^3)$. On comparable nonzero labels it bounds the
difference of the second derivatives by the label increment times
$\abs A^{8/p-3}$, with the corresponding scaled test derivatives.
When the labels are not comparable, use \eqref{res:captest} separately
and split at the size of their increment. Integration of this power
proves the claimed H\"older estimate through zero, since $8/p-2>0$.
No assertion about total variation of derivatives of moving boundaries
is needed.
With $Lf=-\int_Cf\varphi$ and $C_f=\int_C(1-\theta)f\varphi$, put
\begin{equation}\label{res:enlargedcap}
  \mu_A^{\rm lin}=\mu_A+C_S+C_N.
\end{equation}
Use a fixed homogeneous extension of $N$ on $\rho\le2K_0\lambda$. Each added
term has degree eight, for example
$\abs A\int_0^{C\lambda}\rho^{6-m}\dd\rho=C\lambda^8$,
so \eqref{res:captest} remains valid.
The extension of $N$ is made after normalizing its leaf scale. Changing
it changes only a compact test functional with the same degree and
Euler estimates. We use the same extension in every subsequent comparison of the sets, so that its contributions cancel when the exact constant-label
leaf is subtracted.

The vertical moment formula is
\begin{equation}\label{res:moments}
  \int_\R\frac{\abs z^j\,\dd z}{(R^2+z^2)^{(9+\sigma+k)/2}}
  =c_{j,k,\sigma}R^{j-8-\sigma-k},\qquad j<8+\sigma+k.
\end{equation}
For $T=r+\abs t\le Cr^a$, the support radius of $\mu_A^{\rm lin}$
for labels of size $T$ is $o(r)$. Write
\[
K_{x,z}(Y)=(\abs{x-Y}^2+z^2)^{-(9+\sigma)/2}.
\]
Applying \eqref{res:captest} to this test function bounds its second
label difference by
$Cz^2T^{8/p-2}(r^2+z^2)^{-(9+\sigma)/2}$ for $\abs z\le T$.
For $\abs z>T$, total variation and vertical separation give
$CT^{8/p-8-\sigma}\le Cr^{-6-\sigma}T^{8/p-2}$. Thus
\begin{equation}\label{res:pairedcap}
  \kappa_\sigma\left|\int_\R
  [\mu_{A_{\rm out}(t+z)}^{\rm lin}-\mu_{A_{\rm out}(t)}^{\rm lin}]
  (K_{x,z})\,\dd z\right|
  \le C\eps r^{-6-\sigma}T^{8/p-2}.
\end{equation}
Composition with $A_{\rm out}$ preserves the bound: its first derivative
is bounded and $\mu_A'A_{\rm out}''$ has the same or a smaller power.
In the central strip, where no comparison with a constant label is
made, use the absolute estimate
\begin{equation}\label{res:centralcap}
  C\eps\int_\R\frac{\langle s\rangle^{8/p}\,\dd s}
  {(r^2+(s-t)^2)^{(9+\sigma)/2}}
  \le C\eps r^{8/p-8-\sigma},\qquad \abs t\le Cr.
\end{equation}
Both arguments require horizontal separation of the evaluation point
from the support of $\mu_A^{\rm lin}$. Where this separation fails,
including near the collapsed orbit, we use normalized smooth leaf
charts below.

To justify \eqref{res:pairedcap}, note that $K_{x,z}=K_{x,-z}$, so its
integral equals one half of the integral of
\[
[\mu_{A_{\rm out}(t+z)}^{\rm lin}
+\mu_{A_{\rm out}(t-z)}^{\rm lin}
-2\mu_{A_{\rm out}(t)}^{\rm lin}](K_{x,z}).
\]
For labels of size at most $CT$,
\[
\norm{D_Y^kK_{x,z}}_{L^\infty(B_{CT^{1/p}})}
\le Cr^{-k}(r^2+z^2)^{-(9+\sigma)/2},\qquad k\le3.
\]
Every term in \eqref{res:captest} therefore has the same kernel
majorant, since $T^{1/p}/r=o(1)$. The first label derivative cancels
in the displayed symmetric difference. For $\abs z\le T$, integrate
its second derivative using
$\int_\R z^2(r^2+z^2)^{-(9+\sigma)/2}\dd z=Cr^{-6-\sigma}$.
For $\abs z>T$, use the mass bound $C(T+\abs z)^{8/p}$ and vertical
separation. Thus the contribution of the bounded part of the leaf is
estimated directly through $\mu_A^{\rm lin}$, before expanding its
graphical end.

\subsection{Decomposition of the integral}\label{res:phase}

Temporarily use $N(r,A_{\rm out}(t))$ before its cutoff and label correction. Since $p>a$,
$\theta\chi=\chi$ outside a fixed compact set. The height there is
$h=\Psi+\chi W=\ell+N-(1-\chi)W$. Let $E_h$ be the set below this normal height, with
the prescribed bounded completion, and write
$\varpi_h=(\chi_{E_C\times\R}-\chi_{E_h})\,\dd Y$. Then
\begin{align}\label{res:masterphase}
  \varpi_h
  &=\mu_A+L_\theta h+Q_\theta h+\varpi_{\rm comp}\\
  &=L\ell+LN+Q_\theta h+\mu_A^{\rm lin}
  +\int_C(1-\chi)W\varphi\,\dd S+\varpi_{\rm comp}.\nonumber
\end{align}
Here $L_\theta=L+C$, $Q_\theta$ is the tube functional, and
$\varpi_{\rm comp}$ has fixed compact support and finite variation.
The second equality follows from $C_h=C_\Psi=C_S+C_N$.
The noncompact strip of bounded labels is retained, not included in
$\varpi_{\rm comp}$. For singular kernels apply this identity with
smooth truncations, use \eqref{res:pullback} locally, and then remove
cutoffs. Thus the characteristic functions of nonsmooth reference sets
are integrated only against kernels smooth on their support.
More precisely, first put a compact cutoff on the cone coordinates and
regularize the kernel at the evaluation point. The identity is then an
ordinary change of variables between finite integrals. On the local
piece, vary the set $E_h$ and the evaluation point together and subtract
the tangent half-space before expanding. On its complement, the tube
formula and the bounds for $\mu_A$ acting on test functions provide integrable absolute majorants.
They permit removal of the compact cutoff and then the kernel
regularization. This order retains the cancellation of the full
quadratic term and avoids assigning a local curvature to a reference
profile at its label cusp.

The equality $\theta\chi=\chi$ expresses a separation of
scales: where the inner modification is present,
$\rho\gtrsim\abs s^{1/a}$, whereas the support of $\mu_A$ has radius
$O(\abs s^{1/p})$. As $p>a$, these regions are disjoint for large
labels. For bounded labels their failure to be disjoint occurs only in
a fixed compact set, which is exactly the additional distribution
$\varpi_{\rm comp}$. The rest of the bounded-label strip remains in
the displayed tube integrals.

The function $\chi_0(t/r)N(r,A_{\rm eff}(r,t))$ differs from
$N(r,A_{\rm out}(t))$ by its cutoff near $t=0$ and the change of label:
\begin{equation*}
  \norm{(1-\chi_0)N(r,A_{\rm eff})}_{C^{\tau+\alpha}_{\rm sc}}
  \le Cr^{-\beta_N},\qquad \beta_N=6-\mu_N\longrightarrow11/3,
\end{equation*}
\begin{equation*}
  D_A(\rho,s)=N(\rho,A_{\rm eff})-N(\rho,A_{\rm out}(s)),
\end{equation*}
\begin{equation*}
  \norm{D_A}_{C^{\tau+\alpha}_{\rm sc}}\le C\rho^{-\beta_A},\qquad
  \beta_A=7/2<\min(4-\eta,7-\mu_N-(q-d)).
\end{equation*}
For large labels this follows by the mean value formula; on the fixed
completion strip $D_A=\rho^{-6}a(s)+\cdots$, whose scaled
$C^{2,\eta}$ norm is $O(\rho^{-4+\eta})$. Near zero both labels equal
$s$. A normal change of weight $r^{-b}$, $0<b<7$, has curvature weight
$r^{-b-\tau}$ by the local estimate and the radial moment
$\int_0^r\rho^{6-b}\dd\rho$ (or Lemma~\ref{res:coarea} below).
Consequently these changes contribute
\begin{equation}\label{res:fastreturns}
  Cr^{-\beta_N-\tau}+Cr^{-\beta_A-\tau}.
\end{equation}
Here the remote estimate can also be checked directly in cylindrical
coordinates. For $\rho<r/2$, the first variation of the integral is bounded by
\[
C\eps r^{-9-\sigma}
\int_0^r\int_{\abs{s-t}\le Cr}\rho^{6-b}\dd s\,\dd\rho
\le Cr^{-b-1-\sigma}.
\]
For $\rho>2r$, vertical integration leaves
$C\eps\int_r^\infty\rho^{-b-2-\sigma}\dd\rho$.
The remaining vertical tails obey the same moment bound. The local
estimate uses the $C^{\tau+\alpha}$ norms of the two smooth reference
heights, and higher derivatives with respect to these heights have additional small
height-to-radius factors. Thus both comparisons are legitimate global
curvature differences, including the extensions on $\rho\le2K_0\lambda$.

\subsection{Estimates for the change of label and interpolation}\label{res:matchingphase}Split the last integral in \eqref{res:masterphase} with the cutoff
$\theta$. For each vertical coordinate $s$, let
\[
\lambda(s)=(|A_{\rm out}(s)|/a_0)^{1/p},\qquad
\theta_s(\rho)=\theta(\rho/\lambda(s)),
\]
and define the finite signed measure
\[
C_W(s)(\varphi)=\int_C(1-\theta_s(\rho))W(\rho,s)\varphi(P)\,\dd S_P.
\]
Here $W=\ell-S$ with the completed label $A_{\rm out}$; at label zero
the measure is defined by its limit. For large labels its support is
$\rho\le2K_0\lambda(s)$, where $\chi=0$, so it is exactly the part of
$\int_C(1-\chi)W\varphi\,\dd S$ supported in this region. Bounded
labels use the fixed completion in \eqref{res:masterphase}.
Integrating \eqref{res:Wbounds} against $\rho^6\dd\rho$ gives
\begin{equation}\label{res:Wcapmass}
  \norm{C_W(s)}_{\rm TV}\le C\big(\abs s^{q-2}\lambda^{9-n}
  +\abs s^{q_2}\lambda^{7-n}\big)
  \le C\big(\lambda^{10-2p}+\lambda^{7-n+pq_2}\big).
\end{equation}
Define
\begin{equation*}
  \vartheta_W=\max\{(10-2p)/p,(7-n+pq_2)/p\}\longrightarrow4/3.
\end{equation*}
Horizontal separation for nearby labels and vertical separation otherwise give
\begin{equation}\label{res:Wcapsource}
  \left\|\kappa_\sigma\int_\R C_W(s)(K_{x,s-t})\,\dd s
  \right\|_{C^\alpha_{\rm sc}}
  \le C\eps r^{-8-\sigma}T^{\vartheta_W}.
\end{equation}
For large labels the derivatives of the mass bound have powers
$\lambda^{\kappa-jp}$, $j\le2$, for the two exponents $\kappa$ in
\eqref{res:Wcapmass}. The estimate uses the two traces already removed
from $W$. Replacing it by the degree-eight mass bound for
$\mu_A^{\rm lin}$ would lose the required decay. Unlike
\eqref{res:pairedcap}, no second label difference is present in this
integral.

If $\abs{s-t}\le T/2$ and $\abs t\ge Cr$, then $\abs s\asymp T$
and $C_W(s)$ is supported in a ball of radius $O(T^{1/p})=o(r)$.
Its contribution is at most
\[
C\eps T^{\vartheta_W}\int_\R
(r^2+z^2)^{-(9+\sigma)/2}\dd z
=C\eps r^{-8-\sigma}T^{\vartheta_W}.
\]
For $\abs{s-t}>T/2$ and $\abs s\le CT$, use the mass bound
$CT^{\vartheta_W}$ and integrate $\abs{s-t}^{-9-\sigma}$.
For $\abs s>CT$, use $\abs{s-t}\ge c\abs s$ and obtain
\[
C\eps\int_{CT}^\infty s^{\vartheta_W-9-\sigma}\dd s
\le C\eps T^{\vartheta_W-8-\sigma}
\le C\eps r^{-8-\sigma}T^{\vartheta_W}.
\]
The integral converges uniformly since $\vartheta_W<8+\sigma$.
One kernel derivative costs $r^{-1}$; derivatives of the mass for large
comparable labels cost $T^{-1}$. Bounded labels use the finite measure
on a fixed bounded set. The same decomposition therefore proves the
scaled derivative bound, even when the supports for distant labels are
no longer horizontally separated from the evaluation point.
On $\lambda(s)\lesssim\rho\lesssim\abs s^{1/a}$ the end mass is
\begin{equation}\label{res:cutmass}
  C(\langle s\rangle^{\vartheta_1}+\langle s\rangle^{\vartheta_2}),
  \quad \vartheta_1=(9-n)/a+q-2,\quad \vartheta_2=(7-n)/a+q_2.
\end{equation}
Outside the matching ball $D\ge c(r^2+(s-t)^2)^{1/2}$, whence
\begin{equation}\label{res:linearcut}
  \norm{E_{\rm cut,lin}}_{C^\alpha_{\rm sc}}
  \le C\eps r^{-8-\sigma}(T^{\vartheta_1}+T^{\vartheta_2}).
\end{equation}
To see the separation used here, nearby labels have
$\rho\lesssim\abs s^{1/a}$ and hence radial distance comparable to
$r$ at an inner evaluation point outside the matching ball. Where
this radial separation fails, the labels have large vertical distance.
Integrating the two radial monomials in \eqref{res:Wbounds} gives
\eqref{res:cutmass}, and \eqref{res:moments} supplies the remaining
vertical integral. On the matching ball, the local cutoff acts on a
smooth height and is estimated at order $\tau$.
On the outer side $\abs t\ge Cr^a$ the change is $\chi W$, supported
on $\rho\gtrsim\abs t^{1/a}$ for comparable labels. Vertical integration gives
\begin{equation}\label{res:outercut}
  C\eps\big(\abs t^{q-2+(1-n-\sigma)/a}
  +\abs t^{q_2-(n+1+\sigma)/a}+\abs t^{q-n-\tau}\big).
\end{equation}
For example the first term is
$C\eps\abs t^{q-2}\int_{\abs t^{1/a}}^\infty\rho^{-n-\sigma}\dd\rho$;
the last uses homogeneity for $\rho\gtrsim\abs t$. Noncomparable
labels have vertical separation, and cubic terms carry smaller height factors.
This outer estimate is used in place of extrapolating
\eqref{res:linearcut} to the angular scale $\abs t\asymp r^p$.
For $\rho\gtrsim\abs t$, the near-axis expansion of the exact modes
replaces their two-trace expansion. Its homogeneity gives the last
term of \eqref{res:outercut}. For $\rho$ below that scale, the displayed
radial integrals converge at infinity with their stated lower
endpoint. Thus all pieces of the matching change are controlled on
both sides of the interface.

In $Cr\le\abs t\le r^a$, freeze $A=A_{\rm out}(t)$. Exact leaf
minimality and \eqref{res:harmonicell} yield
\begin{align}\label{res:returnidentity}
  \HH_{\rm cyl}(\Psi+\chi W)
  ={}&(\widehat\J_{\rm cyl}-\widehat\J_{\rm rad})N
  -\widehat\J_{\rm cyl}((1-\chi)W)\\
  &+\mathcal Q_{\rm cyl}(\Psi+\chi W)-\mathcal Q_{\rm cyl}(\Psi)
  +\mathcal Q_{\rm cyl}(\Psi)-\mathcal Q_{\rm rad}(\Psi).\nonumber
\end{align}
This is an identity between smooth local integrals and, away from the evaluation point, the identity
\eqref{res:masterphase}, using the functional in \eqref{res:enlargedcap}.
The extra linear term is the full integral of $(1-\chi)W$, whose
parts on $\rho\le2K_0\lambda$ and its complement were estimated in \eqref{res:Wcapsource} and
\eqref{res:linearcut}. Its value at the evaluation point vanishes in
the inner region. In an overlap chart its value at the evaluation point is retained
and is handled by the local matching estimate. Thus the support of
$\chi W$ used in the nonlinear terms is consistent with the cutoff
in identity \eqref{res:masterphase}.
For $(\widehat\J_{\rm cyl}-\widehat\J_{\rm rad})N$, comparable labels contribute
$CT^{\mu_N-2}r^{-6-\sigma}\int_0^r\rho^6\rho^{-6}\dd\rho$;
larger radii give the same power. The contribution of
$\mu_A^{\rm lin}$ is bounded by \eqref{res:pairedcap}.
Thus second label differences and \eqref{res:moments} give
\begin{equation}\label{res:linearreturn}
  \abs{(\widehat\J_{\rm cyl}-\widehat\J_{\rm rad})N}
  \le Cr^{-5-\sigma}T^{\mu_N-2}+C\eps r^{-6-\sigma}T^{8/p-2}.
\end{equation}
For noncomparable labels, $D\ge c(\rho+T+\abs s)$ and the absolute
integral is bounded by
\[
C\eps\int_0^\infty\int_{\abs s\not\asymp T}
\frac{\langle s\rangle^{\mu_N}+T^{\mu_N}}
{(\rho+T+\abs s)^{9+\sigma}}\dd s\,\dd\rho
\le C\eps T^{\mu_N-7-\sigma},
\]
no larger than the first term of \eqref{res:linearreturn}. Derivatives
of the completed label contribute smaller powers.
The local part of
$\mathcal Q_{\rm cyl}(\Psi+\chi W)-\mathcal Q_{\rm cyl}(\Psi)$ is bounded by
\begin{equation}\label{res:staticcubic}
  B(r,T):=Cr^{-3-\sigma-2m}T^2
  (r^{2-n}T^{q-2}+r^{-n}T^{q_2}).
\end{equation}
On smaller radii the term containing only heights at $Y$ gives
\begin{align*}
  I_1&=C\eps r^{-8-\sigma}T^q
  \int_{cT^{1/a}}^{r/2}\rho^{9-n-2p}\frac{\dd\rho}{\rho},\\
  I_2&=C\eps r^{-8-\sigma}T^{q_2+2}
  \int_{cT^{1/a}}^{r/2}\rho^{7-n-2p}\frac{\dd\rho}{\rho}.
\end{align*}
Empty intervals contribute zero. The lower endpoint comes from
$\chi W$; $I_1$ is logarithmic at $\sigma=1$.

\begin{lemma}\label{res:noncomparable}
  The contribution of noncomparable labels to
  $\mathcal Q_{\rm cyl}(\Psi+\chi W)-\mathcal Q_{\rm cyl}(\Psi)$ is bounded in scaled $C^\alpha$ by $C(B+I_1+I_2+r^{-\beta_A-\tau})$.
\end{lemma}
\begin{proof}
  For all powers used here,
  \begin{equation}\label{res:offlabelmoment}
    \int_{\abs s\not\asymp T}\abs s^\mu D^{-9-\sigma-j}\dd s
    \le C(\rho+T)^{\mu-8-\sigma-j},\quad -1<\mu<8+\sigma+j.
  \end{equation}
  Use $H_Y\le C\abs s\rho^{-m}$ and
  $W_Y\le C(\rho^{2-n}\abs s^{q-2}+\rho^{-n}\abs s^{q_2})$.
  For $\abs s\le C\rho$ these value bounds follow from oddness and
  homogeneity: the additional ratios are
  $(\abs s/\rho)^{3-q}$, $(\abs s/\rho)^2$ and
  $(\abs s/\rho)^{2-d}$, all bounded.
  For the term containing only heights at $Y$, on $\rho<r/2$ and
  $\abs s<cT$, integrate up to
  $\min(T,C\rho^a)$. The two radial integrands are
  $\rho^{6-2m-n}\min(T,C\rho^a)^{q+1}$ and
  $\rho^{4-2m-n}\min(T,C\rho^a)^{q_2+3}$, multiplied by
  $C\eps T^{-9-\sigma}$. Below $T^{1/a}$ this gives
  \begin{equation}\label{res:lowlabeltail}
    C\eps(T^{q-8-\sigma+(7-2m-n)/a}
    +T^{q_2-6-\sigma+(5-2m-n)/a}),
  \end{equation}
  with limiting powers $-23/3,-39/5$. Above $T^{1/a}$ it gives $I_1+I_2$.
  For $\abs s>CT$ use \eqref{res:offlabelmoment} and
  $\rho\gtrsim\abs s^{1/a}$; for $\rho\ge r/2$ the resulting limiting
  decays are $23/3$ and $9$.
  All denominators in the integrations below $T^{1/a}$ stay away from
  zero near $\sigma=1$. This is distinct from $I_1$, whose exponent
  approaches zero and whose logarithm is retained. The small-label
  powers in \eqref{res:lowlabeltail}, and the two powers on $\rho\ge r/2$,
  are strictly faster than $r^{-\beta_A-\tau}$ because $T\ge r$.
  For large labels the support condition is important: without
  $\rho\gtrsim\abs s^{1/a}$, the integral would include a
  radial region in which the matching difference is absent.

  For the products in \eqref{res:farcubic} involving heights at $P$, use
  $H_P\le CTr^{-m}$ and $W_P\le C(r^{2-n}T^{q-2}+r^{-n}T^{q_2})$.
  After \eqref{res:offlabelmoment} and radial integration, the terms
  containing $W_Y$ and one or two factors evaluated at $P$ are bounded by
  \[
  C\eps r^{-m}(T^{q-7}+T^{q_2-7}),\qquad
  C\eps r^{-2m}(T^{q-n-\sigma-1}+T^{q_2-n-\sigma-1}).
  \]
  Their ratios to the corresponding terms in $B$ are
  $(r/T)^7$, $(r/T)^9$, $(r/T)^{1+\sigma+n}$ and
  $(r/T)^{3+\sigma+n}$. The radial powers before integration are
  $8-m-n$, $6-m-n$, $9-n$ and $7-n$, all positive.
  Terms containing $W_P$ give
  \[
  C\eps W_P(T^{-2m-1-\sigma}+r^{-m}T^{-m-1-\sigma}
  +r^{-2m}T^{-1-\sigma}),
  \]
  whose ratios to $B$ are
  $(r/T)^{3+\sigma+2m}$, $(r/T)^{3+\sigma+m}$ and $(r/T)^{3+\sigma}$.
  The radial moments converge uniformly. Differentiating partitions costs
  $T^{-1}$, kernels $D^{-1}$ and values at $P$ a factor $r^{-1}$, proving the
  scaled H\"older assertion.
\end{proof}

For comparable labels all remaining products give $B$: their small-radius
moments are $\int_0^r\rho^{5-2m}\dd\rho$,
$\int_0^r\rho^{6-m}\dd\rho$,
$\int_0^r\rho^{7-m-n}\dd\rho$, $\int_0^r\rho^{5-m-n}\dd\rho$,
$\int_0^r\rho^{8-n}\dd\rho$ and $\int_0^r\rho^{6-n}\dd\rho$.
They converge since $m+n=6-\sigma$; large-radius moments include the
kernel decay. The term containing only heights at $P$ is also bounded by $B$.
These moments separate the possible positions of the difference factor
in \eqref{res:farcubic}. A difference at $P$ uses the first
two moments; a difference at $Y$ with one factor evaluated at $P$ uses the next
two, and with two such factors the last two. In particular, the
term containing only values at $P$ is present in this comparison: the heights
at the evaluation point differ by $W_P$, and its bound is
$Cr^{-\sigma-3}H_P^2\abs{W_P}$. This is precisely the size in
\eqref{res:staticcubic}.
On the fixed completion strip the integral containing only heights at $Y$ is
$C\eps\int_0^\infty\rho^{7-3m}(\rho^2+t^2)^{-(9+\sigma)/2}\dd\rho$
$\le C\eps\abs t^{-1-3m-\sigma}$; mixed terms give
$Cr^{-m}\abs t^{-1-2m-\sigma}$ and
$Cr^{-2m}\abs t^{-1-m-\sigma}$, hence are faster.

For $\mathcal Q_{\rm cyl}(\Psi)-\mathcal Q_{\rm rad}(\Psi)$, compare
the label $A_{\rm out}(s)$ with the constant $A_{\rm out}(t)$.
The second label derivative of the term containing only heights at $Y$ is $C\abs A\rho^{-3m-2}$.
The ratios of the $r^{-n}$ term and $N$ to the leading term are
$(\lambda/\rho)^d$ and $(\lambda/\rho)^{6-m}$, also after two label
derivatives. The first vertical moment vanishes; thus
\begin{equation*}
  C\eps r^{-6-\sigma}T\int_{cT^{1/p}}^{r/2}\rho^{4-3m}\dd\rho
  \le C\eps r^{-6-\sigma}T^{8/p-2},
\end{equation*}
Here $m>5/3$, so the small-radius integral is dominated by its lower
endpoint. The identity $1+(5-3m)/p=8/p-2$ converts its power to the
power in \eqref{res:pairedcap}. The part below this lower endpoint
is included in the second label difference of $\mu_A^{\rm lin}$;
the conical integral is used only above that endpoint.
\begin{equation}\label{res:largeradialreturn}
  CT\int_r^\infty\rho^{-3m-2-\sigma}\dd\rho
  \le CTr^{-3m-1-\sigma}.
\end{equation}
One factor evaluated at $P$ gives the latter bound through
$H_P\int_r^\infty\rho^{-2m-2-\sigma}\dd\rho$; on small radii use
$\int_0^r\rho^{5-2m}\dd\rho$.
Two factors evaluated at $P$ give $CT^qr^{-2m-n-1-\sigma}$ or
$CT^{\mu_N}r^{-2m-7-\sigma}$, smaller by
$r^{a(q-1)-d}$ and $r^{a(\mu_N-1)+m-6}$, tending to $r^{-1/6}$ and
$r^{-2/3}$. The terms containing only values at $P$ cancel because
both label choices give the same height at $P$.
The second label derivative of the square of the leading slow height
is constant, which explains the bound with one factor evaluated at $P$.
With two such factors the leading height has zero second label derivative, and
the first nonzero contribution is the $r^{-n}$ term. This is the source
of the extra ratios just displayed. These cancellations occur after
the angular pairing and include the evaluation point; they are not
separate estimates of curvature at a point off a reference surface.
Here $A_{\rm out}'=1+O(T^{q-d-1})$ and
$A_{\rm out}''=O(T^{q-d-2})$. The chain rule adds
$T^{q-d-1}=o(1)$ in the term with one factor evaluated at $P$; with two such factors
the slow branch gives $CT^{q-d}r^{-3m-1-\sigma}$, bounded by
\eqref{res:largeradialreturn} because $q-d<1$. Crossing zero is harmless since
$\abs A^{q-2}$ is locally integrable. Displacement by $W_P$ relative
in the integral of $\mu_A^{\rm lin}$ contributes
\begin{equation*}
  C\eps r^{-9-\sigma}T^{8/p}(r^{2-n}T^{q-2}+r^{-n}T^{q_2}).
\end{equation*}

\subsection{Estimate of the error}\label{res:remainingregions}For $\abs t\le Cr$, use $\abs{\ell/r}\le Cr^{-m}$ and
\eqref{res:harmonicell}. The local cubic bound is $Cr^{-\sigma-3m}$;
the cutoff and label correction of $N$ give \eqref{res:fastreturns},
and $\mu_A^{\rm lin}$ gives \eqref{res:centralcap}.
Remote moments yield the same powers, keeping the exact field $\ell$ intact.
For $r^a\le\abs t\le cr^p$, compare with the exact leaf whose label
is the constant $A_{\rm out}(t)$. The second vertical moment gives
\begin{equation}\label{res:outerbending}
  Cr^{1-n-\sigma}\abs t^{q-2}=Cr^{-2m-\sigma}\abs{t/r^p}^{q-2}.
\end{equation}
To include its tails, scale by $\lambda=(\abs{A_{\rm out}(t)}/a_0)^{1/p}$
and put $R=r/\lambda$, $\delta=\lambda/\abs t$.
The normalized label $\mathfrak a(Z)=A_{\rm out}(t+\lambda Z)/A_{\rm out}(t)$
has second difference $C\delta^2Z^2$ for $\abs Z<c/\delta$.
The leading term is linear in the label. The second moments of the
$r^{-n}$ term, $N$, the nonlinear difference and $\mu_A^{\rm lin}$ are respectively
$CR^{1-n-\sigma}$, $CR^{-5-\sigma}$,
$CR^{-3m-1-\sigma}$ and $C\eps R^{-6-\sigma}$.
They are bounded by $C\langle R\rangle^{1-n-\sigma}$, including the
bounded smooth part of the normalized leaf. Scaling back gives
$\lambda^{-\sigma}\delta^2R^{1-n-\sigma}$, the bound
\eqref{res:outerbending}. The slow relabeling term is smaller by $(r/\abs t)^d$.
For $R\ge K_0$, the three ratios to the second moment of the $r^{-n}$ term are
\[
R^{n-6},\qquad R^{n-3m-2},\qquad \eps R^{n-7}.
\]
Their exponents stay strictly negative near $\sigma=1$. On bounded
normalized charts the second moment is uniformly finite by smoothness.
This justifies the same bound over the whole normalized leaf and not
only on its graphical end.
For $\abs Z>c/\delta$ we estimate the volume between the two surfaces. At a fixed
normalized vertical coordinate $Z$, the label has size at most
$a_Z=C(1+\delta\abs Z)$. For $\rho\ge K_0a_Z^{1/p}$, the volume between the leaf and $C$
in a horizontal ball of radius $\abs Z$ is at most
\[
Ca_Z\int_{ca_Z^{1/p}}^{\abs Z}\rho^{6-m}\dd\rho
\le Ca_Z\abs Z^{7-m},
\]
with an empty interval contributing zero. The part inside $\rho\le K_0a_Z^{1/p}$ has volume
$Ca_Z^{8/p}$. If its radius exceeds the ball radius, this second
bound alone controls the whole intersection. On the large-$Z$ range
in question the observation point satisfies
$R/\abs Z\le C r/\abs t=o(1)$, so its horizontal displacement does
not change these kernel bounds.

In each horizontal shell $\rho\ge\abs Z$, the volume between the leaf and $C$
is $Ca_Z\rho^{7-m}$ and the kernel is bounded by
$C\rho^{-9-\sigma}$. The dyadic shell sum is therefore
$Ca_Z\abs Z^{-2-m-\sigma}$, which is the same bound obtained from
the inner horizontal ball. Finally, the substitution $u=\delta Z$
makes the parameter dependence explicit:
\begin{align*}
  \int_{c/\delta}^\infty(1+\delta Z)Z^{-2-m-\sigma}\dd Z
  &=\delta^{p+\sigma}\int_c^\infty(1+u)u^{-2-m-\sigma}\dd u,\\
  \int_{c/\delta}^\infty(1+\delta Z)^{8/p}Z^{-9-\sigma}\dd Z
  &=\delta^{8+\sigma}\int_c^\infty(1+u)^{8/p}u^{-9-\sigma}\dd u.
\end{align*}
Both integrals on the right are uniformly finite. Including the
normalization, the large-label tail is consequently bounded by
$C\eps(\delta^{p+\sigma}+\delta^{8+\sigma})$, with its full
horizontal contribution.
Here $p+\sigma>2$, and in conical coordinates the absolute bound is
$C\eps\abs t^{-m-\sigma}$. At $\abs t\asymp r^p$, smooth charts on the whole leaf give $Cr^{-\sigma}r^{2-2p}=Cr^{-2m-\sigma}$.
Thus Taylor expansion is used only where the normalized label varies
by a bounded relative amount. Outside that range the volume estimate is smaller than the second-moment scale. On the collapsed
orbit the normalized geometry is smooth, so the same argument applies
without conical coordinates; the label changes by $O(r^{1-p})$
across a chart of radius $r$.
The inner modification is \eqref{res:outercut}.

An ambient annulus of radius $D\ge cr^p$ contains a region between the approximate graph and $\Sigma$ of volume at most
$CD\int_0^{CD}\rho^7\min(1,D/\rho^p)\dd\rho\le CD^{10-p}$.
Dyadic kernel integration therefore gives
\begin{equation*}
  C\eps r^{-p(p+\sigma-1)}.
\end{equation*}
The compact completion adds $C\eps r^{-9-\sigma}$.

\begin{lemma}\label{res:holder}
  All these estimates hold in scaled $C^\alpha$, uniformly as $\sigma\to1$.
\end{lemma}
\begin{proof}
  Locally use Section~\ref{sec:preliminaries}, with the tangent-half-space
  localization of \cite[Section 8]{DDPW}. The remote bounds are verified
  in the present geometry as follows. Localize each remote region by a smooth
  partition $\chi_j(X,Y)$ with $\abs{D_X\chi_j}\le C/r$. A regional
  bound $M_j$ and a differentiated bound $CM_j/r$ imply a scaled
  $C^\alpha$ bound $CM_j$ by interpolation. Both terms in the derivative
  of $\chi_j$ times the kernel are included. On comparable labels,
  each additional label derivative costs $T^{-1}\le r^{-1}$; on the
  remaining labels one differentiates the regular test kernel instead.
  The Euler estimates apply to this differentiated test function and
  the integrated third label estimate supplies its H\"older remainder.
  For $N$ use the same argument with \eqref{res:Nderivatives}.
  The absolute dyadic sums above converge with these derivative powers,
  so smooth regularizations can be removed in the asserted norm.
  In $I_1$ absorb the logarithm by a fixed arbitrarily small loss of a power.
\end{proof}

\begin{lemma}\label{res:coarea}
  For the approximate graph and its small normal perturbations,
  \begin{equation*}
    \int_{\mathcal S\cap B_D(X)}r(Y)^{-b}\dd S_Y\le C_bD^{8-b},
    \qquad D\ge r(X)/2,\quad 0\le b<7.
  \end{equation*}
  A normal change of weight $r^{-b}$ therefore has remote curvature
  weight $r^{-b-\tau}$, also for differentiated kernels.
\end{lemma}
\begin{proof}
  In the angular region the homogeneous leaf label has the form
  $\rho^pg(\vartheta)$; the vertical height is obtained by the monotone
  change $A_{\rm out}^{-1}$. The classical profile is monotone away from
  its collapsed orbit and has a nondegenerate extremum there
  \cite[Lemmas 2.1 and 9.2]{DPKW}. The smooth normalized fractional leaf
  preserves these properties. Thus its orbital angular measure has
  bounded pushforward by $g$: in the interior $g'\ne0$, and at the axis
  $g(0)-g(\vartheta)\asymp\vartheta^2$ while the orbital density is
  $O(\vartheta^3)$. At the cone use the differentiated leaf expansion.
  The bounded derivative of the change of label and the matching bounds
  therefore give angular measure $C\min(1,D/\rho^p)$ for a vertical
  interval of length $2D$. Consequently
  \[
  \int_{\mathcal S\cap B_D(X)}r^{-b}\dd S
  \le C\int_0^{CD}\rho^{p+6-b}\min(1,D/\rho^p)\dd\rho
  \le C_bD^{8-b},
  \]
  splitting at $D^{1/p}$. Cylindrical charts and small normal changes
  preserve this bound. More explicitly, the lower integral has size
  $CD^{(p+7-b)/p}$ and the upper one is at most
  $CD\int^{CD}\rho^{6-b}\dd\rho=CD^{8-b}$. The lower power is no
  larger because $p>1$ and $b<7$. On the compact part replace $r$ by a
  smooth positive weight. The region between the normal graphs has volume $CD^{8-b}$ on a remote
  annulus; multiplying by $\eps D^{-9-\sigma}$ and summing proves the consequence.
\end{proof}

\begin{proposition}\label{res:assembled}The surface $\cP_\sigma$ satisfies, for $r\ge R_0$ and $\sigma$ close to one,
  \begin{equation}\label{res:principalbound}
    \norm{\HH_\sigma(\cP_\sigma)}_{C^\alpha_{\rm sc}(\mathcal U_r)}
    \le Cr^{-\gamma-\tau}.
  \end{equation}
  Its scaled charts have uniform bounds for every fixed finite number of
  derivatives, and, with the upward orientation,
  \begin{equation}\label{res:geometry}
    \abs{A_{\cP_\sigma}}\le C/r,\qquad
    cr^{-m}\le(\nu_{\cP_\sigma})_9\le Cr^{-m}.
  \end{equation}
\end{proposition}
\begin{proof}
  In the intermediate region the complete bound is the sum
  \begin{align*}
    C\big[&r^{-5-\sigma}T^{\mu_N-2}+B+Tr^{-3m-1-\sigma}
    +\eps r^{-6-\sigma}T^{8/p-2}\\
    &+\eps r^{-8-\sigma}\!\sum_{j=W,1,2}T^{\vartheta_j}+I_1+I_2
    +r^{-\beta_N-\tau}+r^{-\beta_A-\tau}+\eps r^{-p(p+\sigma-1)}\\
    &+\eps r^{-9-\sigma}T^{8/p}(r^{2-n}T^{q-2}+r^{-n}T^{q_2})
    +\eps r^{-9-\sigma}\big].
  \end{align*}
  This follows from \eqref{res:masterphase} and
  \eqref{res:returnidentity}: $\widehat\J_\Sigma^\sigma\ell=0$,
  the quadratic variation vanishes, $\mu_A^{\rm lin}$ is estimated by
  its second label difference, and the integral of $C_W(s)$ is counted
  once. The estimates above include every term containing heights at the
  evaluation point. The following
  table records the least decay powers at $\sigma=1$; intermediate powers
  are maximized over $r\le T\le Cr^{5/2}$.
  \begin{center}
    \begin{tabular}{@{}p{0.69\linewidth}l@{}}
      Contribution & Decay power\\ \hline
      $(\widehat\J_{\rm cyl}-\widehat\J_{\rm rad})N$ on the graphical part
      & $31/6$\\
      Local integral of $\mathcal Q_{\rm cyl}(\Psi+\chi W)-\mathcal Q_{\rm cyl}(\Psi)$
      & $17/3,6$\\
      $\mathcal Q_{\rm cyl}(\Psi)-\mathcal Q_{\rm rad}(\Psi)$ on the graphical part
      & $11/2,16/3$\\
      Second label difference of $\mu_A^{\rm lin}$ & $16/3$\\
      Integral of $C_W(s)$ & $17/3$\\
      $\int_C\theta(1-\chi)W\varphi\,\dd S$ & $14/3,5$\\
      $I_1,I_2$ & $17/3$ (logarithm), $6$\\
      Cutoff and change of label in $N$ & $17/3,11/2$\\
      Change of evaluation point in the integral of $\mu_A^{\rm lin}$ & $6,19/3$\\
      Ambient annuli $|X-Y|\ge cr^p$ & $9$\\
      Region $|t|\le Cr$ & $17/3,7,19/3$\\
      Regions $|t|\asymp r^a$ and $r^a\le|t|\le cr^p$ & $14/3$\\
      Region $|t|\asymp r^p$, including the collapsed orbit & $5$
    \end{tabular}
  \end{center}
  The estimates for noncomparable labels in Lemma~\ref{res:noncomparable}
  are absorbed by the displayed bounds for $B$, $I_1$, $I_2$ and
  $r^{-\beta_A-\tau}$.
  For the overlap, the two matching height exponents are
  $n-2+a(2-q)\to8/3$ and $n-aq_2\to3$; a local cutoff costs
  $r^{-\tau}$. Formula~\eqref{res:outercut} has powers
  $\abs t^{-28/15},\abs t^{-2},\abs t^{-11/3}$ at $\sigma=1$ and
  is largest at $\abs t=r^a$. These give the stated outer bounds.
  The smallest power is $14/3$, with strict margin
  \[
  \frac{14}{3}-(\gamma+2)=\frac1{24}.
  \]
  The logarithmic loss occurs in a faster term. Continuity in $m,n,\sigma$
  preserves all inequalities uniformly, and Lemma~\ref{res:holder} proves
  \eqref{res:principalbound}.

  Finally, the inner height divided by $r$ tends to zero with differentiated
  bounds. The term $N(r,A_{\rm eff})$ is used only for $\abs t\ge r$, where each label
  derivative costs at most $C/r$; the exact modes are smooth through zero.
  The overlap follows from \eqref{res:Wbounds}, and the angular region
  from smooth charts on the normalized leaf. This proves all finite scaled
  geometry bounds. The positive label derivative and its leading term
  $r^{-m}$ give the normal estimate in \eqref{res:geometry}; the same
  charts give the curvature bound. The fixed compact completion uses a
  smooth positive radius comparable to one.
\end{proof}

\section{The classical core and compact matching}\label{glu:coresection}

Let $F_{\rm BDG}$ be the entire minimal graph given by
\cite[Lemma~2.1 and Theorem~2]{DPKW}, following \cite{BDG}. We use
its rotational symmetry, its oddness under $u\leftrightarrow v$, and
the comparison
\[
F_{\rm BDG}(r,\theta)\asymp F_0(r,\theta)=r^3g(\theta)
\quad\hbox{in }\{v>u\},\qquad r\gg1.
\]
Here $F_0$ is the function in \eqref{eq:classicalprofile}.
Our classical foliation is fixed by $M_0$ in
\eqref{mod:classicalleaf}. Its homogeneous label $F_{\infty,1}$ is
comparable to $F_0$: both angular functions are positive in the
sector and vanish linearly at the cone. This comparison suffices
below. We construct a classical minimal graph whose normal difference
from the prescribed end $\mathcal P_1$ has the decay needed for
compact matching.

\begin{proposition}[Compatibility of the classical core]\label{glu:core}
  There is a smooth entire classical minimal graph $\Gamma_*$, with the
  symmetries of $\mathcal P_1$, whose end is a single normal graph over
  $\mathcal P_1$ with height $R_1$ satisfying
  \begin{equation}\label{glu:fastcore}
    \norm{R_1}_{X_{1,\beta_1}(r\ge R_*)}<\infty,
    \qquad \beta_1=\frac{31}{12}.
  \end{equation}
  Moreover, $(\nu_{\Gamma_*})_9\asymp r^{-2}$ and its tangent charts
  have the uniform scaled bounds of Proposition~\ref{res:assembled}.
\end{proposition}

\begin{proof}
  All normals in this proof point upward. Write $F_P$ for the graph
  function of $\mathcal P_1$ and $w=(\nu_{\mathcal P_1})_9$. On the end,
  Proposition~\ref{res:assembled} gives
  \begin{equation}\label{glu:classicalgeometry}
    w\asymp r^{-2},\qquad |A|\le C/r,\qquad
    \norm{\mathcal H_1(\mathcal P_1)}_{C^\alpha_{\rm sc}}
    \le Cr^{-\gamma-2},\quad \gamma=\frac{21}{8}.
  \end{equation}
  Let $B=B_{{\rm end},1}$ be the positive exterior barrier constructed
  in Section~\ref{glu:endbarriers}, with weight $b=\beta_1$. Its
  construction uses only the principal end and the classical smoothing,
  so is independent of the core. Thus
  \begin{equation*}
    B\asymp r^{-b},\qquad J_{\mathcal P_1}B\le-cr^{-b-2},
    \qquad |\nabla^jB|\le C_jr^{-b-j}.
  \end{equation*}
  For $C\ge C_0$ and $r\ge R$, the normal graphs
  $S_C^\pm=\mathcal P_1[\pm CB]$ are strict barriers with opposite
  curvature signs, provided $CR^{-b-1}$ is small and $R$ is large.
  Indeed, their linear term has size $Cr^{-b-2}$, their quadratic error
  is at most $C' C^2r^{-2b-3}$, and the residual in
  \eqref{glu:classicalgeometry} is smaller since $b<\gamma$.

  These barriers are vertical graphs when $CR^{1-b}$ is sufficiently
  small, because their horizontal projection satisfies
  \begin{equation}\label{glu:fixedbarriergraph}
    |D_x(CB\nu_{\rm hor})|
    \le C' w^{-1}(|\nabla(CB)|+CB|A|)
    \le C'C r^{1-b}.
  \end{equation}
  Their graph functions satisfy
  \begin{equation}\label{glu:verticalnormalvariation}
    F_{S_C^+}-F_P\asymp Cr^{2-b},\qquad
    F_P-F_{S_C^-}\asymp Cr^{2-b}.
  \end{equation}
  To verify this, parametrize the normal deformation by
  $X+sCB(X)\nu(X)$, $0\le s\le1$, and write $F_s$ for its vertical
  graph function. Its horizontal projection is
  $T_s(x)=x+sCB(x)\nu_{\rm hor}(x)$.
  At a fixed horizontal point $y$, set $q_s(y)=T_s^{-1}(y)$.
  If $\nu_s$ is the upward normal of the deformed graph and
  $w_s=(\nu_s)_9$, differentiation at fixed $y$ gives
  \[
  \partial_sF_s(y)
  =\frac{CB(q_s(y))\,\nu(q_s(y))\cdot\nu_s(y)}{w_s(y)}.
  \]
  At $s=0$ this gives $\partial_sF_s=CB/w$, the relation between
  normal and vertical perturbations in \cite[(2.10)]{DPKW}. Estimate \eqref{glu:fixedbarriergraph} keeps $q_s(y)$ at
  a comparable horizontal radius. The normal change is
  $O(Cr^{-b-1})=o(r^{-2})$, so $w_s\asymp w\asymp r^{-2}$ and
  the scalar product of the normals is bounded away from zero.
  Integration in $s$ proves \eqref{glu:verticalnormalvariation}.
  The normal tube and horizontal inverse are uniform on a slightly
  reduced exterior region by the scaled geometry. We always define the
  barriers on a slightly larger annulus than the comparison domain.

  For $L\to\infty$, solve the minimal graph Dirichlet problem in the
  horizontal ball $B_L\subset\R^8$ with boundary value $F_P$.
  Existence and uniqueness follow from \cite{JS}; denote the solution by
  $F_L$. Its symmetries are preserved. Fixed ambient dilations of $F_{\rm BDG}$ bracket the boundary
  values in the sector $T=\{v>u\}$ for every large $L$. We verify the
  comparison with $F_0$. In the inner and matching regions,
  \[
  \psi(r,t)=tr^{-2}(1+o(1))
  \]
  uniformly: the central cylindrical terms have relative size
  $O(r^{-2/3})$, and $r^{-3}\phi_{4/3}(t)$ has relative size
  $O(|t|^{1/3}/r)=O(r^{-1/6})$. The conical normal coordinate divided
  by $r$ is comparable to $\theta-\pi/4$. Since
  $g(\theta)\asymp\theta-\pi/4$ near the cone, we obtain
  $F_P\asymp r^3g(\theta)$ there. In the outer region the classical
  leaf parameter is $F_{\infty,1}$, and
  \[
  A_{\rm out}(F_P(r,\theta))=F_{\infty,1}(r,\theta)\asymp F_0(r,\theta).
  \]
  The relation $A_{\rm out}(t)=t+o(t)$ gives the same comparison on the
  remaining angular intervals, including the smooth part at the axis.
  Thus $F_P\asymp F_0$. Ambient dilation of $F_{\rm BDG}$ multiplies
  its leading cubic coefficient by the inverse square of the dilation
  factor. Choosing one small and one large factor gives the required
  upper and lower boundary inequalities. Both comparison graphs vanish
  on $\partial T$. Comparison gives
  bounds for $F_L$ on every fixed ball, uniformly in $L$. Interior minimal
  graph estimates therefore give smooth compact convergence along a
  subsequence to an entire minimal graph $F_*$.

  We obtain the stronger end estimate before taking this limit. Fix the
  strict barriers $S_{C_0}^\pm$ outside a fixed radius $R_0$. The compact
  bounds just obtained allow vertical translates $S_{C_0}^++M$ and
  $S_{C_0}^--M$ to bracket every $F_L$ at the inner boundary. They also
  bracket the data $F_P$ at $r=L$ and zero on the conical boundary.
  Finite-annulus comparison consequently gives a constant $M'$ independent
  of $L$ such that
  \begin{equation*}
    |F_L-F_P|\le M'\qquad(R_0\le r\le L).
  \end{equation*}
  Choose now a large fixed $R$ and
  $C_R=K(M'+1)R^{b-2}$, with fixed $K$ large. By
  \eqref{glu:verticalnormalvariation}, the unshifted barriers
  $S_{C_R}^\pm$ bracket $F_L$ at $r=R$. They remain strict vertical
  barriers, since
  \[
  C_RR^{1-b}=K(M'+1)R^{-1}\longrightarrow0,
  \qquad C_RR^{-b-1}=K(M'+1)R^{-3}\longrightarrow0.
  \]
  At the outer boundary they bracket $F_L=F_P$, and on the symmetry
  boundary their values have the required opposite signs. Comparison on
  $R<r<L$, followed by $L\to\infty$, traps $\Gamma_*={\rm graph}(F_*)$
  between $S_{C_R}^-$ and $S_{C_R}^+$. In particular, the barriers $S_{C_R}^\pm$
  are compared on bounded domains with prescribed outer data; no
  comparison principle at infinity is required.

  Finally this trapping gives a single normal correspondence on the end.
  The principal surface has a normal tube of radius $cr$. Indeed, the
  inner region is a small scaled graph over $C\times\R$, whose orbital
  radii $(r-s)/\sqrt2,(r+s)/\sqrt2$ determine its normal coordinates
  uniquely. The normalized whole leaf has positive reach on compact
  pieces and a conical normal tube on its end. Its scale changes by
  $O(r/|t|)$ across an outer chart. These descriptions agree on the
  overlap; points at distance less than $cr$ have comparable horizontal
  radii and vertical difference less than $cr$, excluding a second
  distant piece from the tube.

  The trapping width is $O(r^{-b})$. In a ball of radius $c_0r$ the
  principal surface lies in a tangent slab of width $Cc_0^2r$, and the
  trapped minimal graph lies in a slab of width $Cc_0^2r+Cr^{-b}$.
  Its normalized flatness is $Cc_0+Cr^{-b-1}/c_0$. Choose first $c_0$
  small and then $r$ large. Minimal graphs are locally area minimizing,
  so classical improvement of flatness \cite{DeGiorgi} makes $\Gamma_*$
  a single sheet transverse to the principal normal fibers.
  Each such fiber meets the graph because the two trapping boundaries
  lie on opposite sides of the graph, and transversality gives uniqueness. Hence its
  normal height obeys $|R_1|\le Cr^{-b}$. In a scaled chart the two
  surfaces are uniformly close graphs, so subtract their classical
  minimal surface equations. The difference satisfies a uniformly
  elliptic equation obtained by averaging the linearization along the
  segment joining the two graph functions. Its coefficients have uniform
  local regularity, and its forcing is the curvature of $\mathcal P_1$.
  Interior estimates on a smaller chart give
  \[
  r^b\|R_1\|_{C^{2,\alpha}_{\rm sc}}
  \le C\big(r^b\|R_1\|_{L^\infty}
  +r^{b+2}\|\HH_1(\mathcal P_1)\|_{C^\alpha_{\rm sc}}\big)
  \le C(1+r^{b-\gamma}).
  \]
  The last expression is bounded because $b<\gamma$. Differentiating
  the same equation on successively smaller charts gives the fixed
  higher derivatives used in the compact change of variables. This proves the
  full estimate \eqref{glu:fastcore}, rather than just a bound on the
  distance between the two surfaces. In particular,
  $|\nu_{\Gamma_*}-\nu_{\mathcal P_1}|\le Cr^{-b-1}=o(r^{-2})$,
  which preserves the vertical-normal comparison and the scaled geometry.
\end{proof}

\subsection{Compact matching}\label{glu:diagonalsection}

\begin{lemma}[Diagonal compact matching]\label{glu:diagonal}
  There are smooth entire graphs $\Sigma_\sigma$, equal to
  $\mathcal P_\sigma$ outside a radius tending to infinity, such that
  \begin{equation*}
    \norm{\mathcal H_\sigma(\Sigma_\sigma)}_{Y_{\sigma,\beta}}
    \longrightarrow0,
  \end{equation*}
  and their scaled geometry is uniform. They converge to $\Gamma_*$ on
  compact sets and satisfy
  \begin{equation*}
    cr^{-m_\sigma}\le(\nu_{\Sigma_\sigma})_9\le Cr^{-m_\sigma}.
  \end{equation*}
  Any finite collection of compact continuity requirements for the
  geometry and the barriers can be imposed in the same choice.
\end{lemma}

\begin{proof}
  At $\sigma=1$, insert the core by cutting off its normal remainder:
  $\Sigma_{1,R}$ equals $\Gamma_*$ on $r\le R$, has height
  $\chi_RR_1$ over $\mathcal P_1$ on the end, and equals
  $\mathcal P_1$ on $r\ge2R$. The transition is a vertical graph by
  \eqref{glu:fixedbarriergraph}. The cutoff commutator is
  $O(R^{-\beta_1-2})$ and the outer principal residual is
  $O(r^{-\gamma-2})$. Thus
  \begin{equation*}
    \norm{\mathcal H_1(\Sigma_{1,R})}_{Y_{1,\beta}}
    \le C(R^{\beta-\beta_1}+R^{\beta-\gamma})
    \le CR^{-1/12}.
  \end{equation*}
  For fixed $R$, choose a smooth change of variables which identifies
  $\mathcal P_1$ with $\mathcal P_\sigma$ on $r\le3R$ and converges
  to the identity there. Apply it to $\Sigma_{1,R}$ on $r\le2R$,
  and keep $\mathcal P_\sigma$ on $r\ge2R$. The two definitions agree
  in the overlap, where $\Sigma_{1,R}=\mathcal P_1$. Denote the resulting
  surface by $\Sigma_{\sigma,R}$. Normalized fractional curvature converges
  to classical curvature in compact $C^\alpha$ by
  \cite[Appendix~A]{DDPW}. The modified subgraphs differ only in a bounded ambient
  set of volume $V_R$; at remote evaluation points their contribution is
  bounded in scaled $C^\alpha$ by
  \begin{equation*}
    C(1-\sigma)V_Rr^{-9-\sigma},\qquad r\ge4R.
  \end{equation*}
  For completeness, multiplication of this tail by the source weight gives
  \[
  r^{\beta+\tau}(1-\sigma)V_Rr^{-9-\sigma}
  =(1-\sigma)V_Rr^{\beta-8}.
  \]
  Since $\beta<8$, its supremum on $r\ge4R$ tends to zero for every
  fixed $R$. A scaled derivative of the kernel gives the same bound,
  and hence also the scaled H\"older seminorm. On the unchanged end,
  Proposition~\ref{res:assembled} contributes at most
  $CR^{\beta-\gamma}=CR^{-1/8}$. In the transition region, the
  classical commutator has weight $R^{\beta-\beta_1}=R^{-1/12}$,
  while every coefficient on this fixed region varies continuously with
  $\sigma$. On the inner compact region the limiting surface is exactly
  minimal. Taking the maximum of these three bounds proves that the
  resulting surfaces satisfy
  \begin{equation*}
    \norm{\mathcal H_\sigma(\Sigma_{\sigma,R})}_{Y_{\sigma,\beta}}
    \le CR^{-1/12}+\omega_R(1-\sigma),\qquad
    \omega_R(\varepsilon)\longrightarrow0
    \quad\hbox{for each fixed }R.
  \end{equation*}
  Choose $R_j\to\infty$ and then $\varepsilon_j\downarrow0$ so that
  $\omega_{R_j}(\varepsilon)\le1/j$ for
  $0<\varepsilon\le\varepsilon_j$. Decrease each $\varepsilon_j$ further
  so that, on $r\le4R_j$, the geometry and the corresponding barrier retain
  the required strict inequalities and
  $r^{m_\sigma}(\nu_{\Sigma_{\sigma,R_j}})_9$ stays between fixed
  positive bounds. These are finitely many continuity requirements on a
  fixed compact set; outside it the uniform principal estimates apply.
  Taking $R(\varepsilon)=R_j$ for
  $\varepsilon_{j+1}<\varepsilon\le\varepsilon_j$ proves all assertions.
  No rate of compact convergence is needed in this choice. The powers
  $R^{-1/12}$ and $R^{-1/8}$ are obtained before the parameter is
  varied; each remaining error is then made small on a fixed compact
  region. In particular, no classical asymptotic expansion is being
  continued to infinity at a different fractional exponent.
\end{proof}

\section{The inverse and the nonlinear correction}\label{glu:inverse}

Write $\Sigma_\sigma=\Sigma_{\sigma,R(\sigma)}$ for the completed
approximation and $L_\sigma=L_{\Sigma_\sigma,\sigma}$. Thus
\begin{equation*}
  L_\sigma u(X)=\PV\int_{\Sigma_\sigma}(u(Y)-u(X))K_\sigma(X,Y)\dd S_Y
  +V_\sigma(X)u(X),
\end{equation*}
where $K_\sigma(X,Y)=\kappa_\sigma|X-Y|^{-8-\tau}$ is positive and
symmetric. The potential is
\[
V_\sigma(X)=\int_{\Sigma_\sigma}
(1-\nu(X)\cdot\nu(Y))K_\sigma(X,Y)\dd S_Y.
\]
In a tangent ball of radius $cr$, its numerator is $O(|X-Y|^2/r^2)$.
Integration against the normalized kernel, followed by the remote
estimate below with weight zero, gives
\begin{equation}\label{glu:potentialbound}
  0\le V_\sigma(X)\le Cr^{-\tau},\qquad
  \|r^\tau V_\sigma\|_{C^\alpha_{\rm sc}}\le C.
\end{equation}
The H\"older assertion follows by the same local cancellation and by
differentiating the remote kernel. The approximation has nonzero
fractional curvature; we therefore construct a strict positive barrier
for this operator. The vertical Jacobi field is used on the exact
classical core in the intermediate construction below.

\subsection{Positive barriers}\label{glu:endbarriers}

We construct the barriers for $b\in\{b_s,\beta,\beta_1\}$.
In the region near $C\times\R$ take $B_{\rm in}=r^{-b}$. Its image
under the cylindrical Jacobi operator is
$I_\sigma(b)r^{-b-\tau}\le-cr^{-b-\tau}$.

In the outer region use the strict barrier on the exact fractional
leaf from \eqref{mod:exactleafbarrier}. The normal change of variables and its error were estimated after the
construction of $M_\sigma$. The horizontal section at height $t$ is
the leaf with parameter $A_{\rm out}(t)$. Its scale is
$\lambda(t)=(|A_{\rm out}(t)|/a_0)^{1/p}\asymp|t|^{1/p}$. Define
\begin{equation}\label{glu:outerbarrier}
  B_{\rm out}(r,t)=c_b^{-1}\lambda(t)^{-b}
  B_{{\rm leaf},\sigma}(r/\lambda(t)).
\end{equation}
Scaling gives its strict source on the leaf with its label $A$ fixed. In the overlap
$|t|\asymp r^a$, $a=5/2$, its leading term equals $r^{-b}$ and
\begin{equation}\label{glu:barriermismatch}
  B_{\rm out}-B_{\rm in}=O_{C^k_{\rm sc}}(r^{-b-\eta_b}),
  \qquad \eta_b=(p-a)(3-b)/p>0.
\end{equation}
Indeed the remainder in \eqref{glu:outerbarrier} is
$O(\lambda^{3-b}r^{-3})$, and $\lambda\asymp r^{a/p}$ there.
Join the two positive barriers with the geometric partition. On a chart
of radius $cr$, the deviation from the cylinder is $O(r^{a-p})$,
and the relative change of the outer leaf scale is $O(r/|t|)$.
Comparison of kernels and normal potentials, and the cutoff
commutator estimated by \eqref{glu:barriermismatch}, therefore give
an error bounded by
\[
C(r^{a-p}+r^{1-a}+r^{-\eta_b})r^{-b-\tau}.
\]

The complement of the tangent chart is estimated directly on $\Sigma_\sigma$. The coarea geometry of Lemma~\ref{res:coarea} gives
\begin{equation}\label{glu:coareaslice}
  \HH^7(\Sigma_\sigma\cap\{t=s,\ \rho<|x|<2\rho\})\le C\rho^7,
  \qquad |\nabla_{\Sigma_\sigma}t|\ge c
\end{equation}
on the end. For a source annulus $\rho<r(X)/2$, integration first in
the vertical variable gives
\[
\int_\R(r^2+(t-s)^2)^{-(9+\sigma)/2}\dd s
=D_\sigma r^{-8-\sigma}.
\]
The inner radial annuli consequently contribute
\begin{equation}\label{glu:innerbarriertail}
  C\eps r^{-8-\sigma}\sum_{\rho<r/2}\rho^{7-b}
  \le C\eps r^{-b-\tau}.
\end{equation}
The sum is dyadic and its lower endpoint is harmless precisely when
$b<7$. For $\rho>2r$, the corresponding sum is
$C\eps\sum_{\rho>2r}\rho^{-b-\tau}$, which converges because
$b>-\tau$. Annuli of comparable horizontal radius outside the tangent
ball are bounded using \eqref{glu:coareaslice}. Evaluation-point
terms have the same bounds with $B(X)$ outside the integral. The
potential tail is the weight-zero calculation multiplied by $B(X)$,
since $0\le1-\nu(X)\cdot\nu(Y)\le2$. Fixed compact pieces contribute $O(\eps r^{-9-\sigma})$. Thus the remote error is $C\eps r^{-b-\tau}$.
This estimate includes distant parts of the two regions without identifying
them pointwise with a model surface.

We have constructed
\begin{equation*}
  \begin{split}
    B_{{\rm end},\sigma}&\asymp r^{-b},\\
    L_\sigma B_{{\rm end},\sigma}
    &\le-[c-C(r^{a-p}+r^{1-a}+r^{-\eta_b}+\eps)]r^{-b-\tau}.
  \end{split}
\end{equation*}
First enlarge the fixed end threshold, and then decrease $\eps$, so
that the bracket is positive. At $\sigma=1$, $b=\beta_1$, this
construction uses only the classical leaf and the principal end, and
supplies the barrier used in Proposition~\ref{glu:core}.

\label{glu:classicalbarrier}We next complete the barrier across the bounded core. On the exact
classical graph $\Gamma_*$ the positive function
$w_*=(\nu_{\Gamma_*})_9$ satisfies $J_{\Gamma_*}w_*=0$.
Transfer the exterior classical barrier to $\Gamma_*$ by its
normal correspondence with $\mathcal P_1$, whose height belongs to
$X_{1,\beta_1}$. Its strict inequality
persists by \eqref{glu:fastcore}. Put
$U=B_{{\rm end},1}/w_*\asymp r^{2-b}$, extend it positively inside,
and choose a small regular value $\delta>0$ such that
$D=\{U>\delta\}$ is compact and its boundary lies on the end.
The outward conormal derivative $\partial_\nu U$ is strictly negative
on $\partial D$.

Solve the uniformly elliptic problem on this fixed domain
\begin{equation}\label{glu:corebarrierproblem}
  -\operatorname{div}_{\Gamma_*}(w_*^2\nabla u)=w_*r^{-b-2},
  \qquad u=0\quad\hbox{on }\partial D.
\end{equation}
The maximum principle gives $u>0$, and the identity
$J(w_*u)=w_*^{-1}\operatorname{div}(w_*^2\nabla u)$ gives
$J(w_*u)=-r^{-b-2}$. The functions
$B_c=w_*(K\delta+u)$ and $B_e=KB_{{\rm end},1}$ agree on the
boundary. For large fixed $K$ their outward derivatives satisfy
\begin{equation*}
  \partial_\nu B_e-\partial_\nu B_c
  =w_*(K\partial_\nu U-\partial_\nu u)<0.
\end{equation*}
Extend $B_c$ across a collar, preserving its value, its derivatives through order two at the
boundary, and its strict inequality. The derivative separation makes the minimum
select $B_c$ on the inner side and $B_e$ on the outer side.

To smooth this minimum, choose an increasing concave smooth function
$H(s_1,s_2)$, homogeneous of degree one and equal to the smaller
argument outside a narrow ratio interval about one. Concavity gives
\[
H(B_1(Y),B_2(Y))-H(B_1(X),B_2(X))
\le\sum_{i=1}^2H_i(B_1(X),B_2(X))(B_i(Y)-B_i(X)).
\]
After integration, Euler's identity $\sum H_iB_i=H$ includes the
potential term and yields
\begin{equation}\label{glu:concavemin}
  L H(B_1,B_2)\le H_1LB_1+H_2LB_2.
\end{equation}
For the classical operator the same inequality follows from the
nonpositive Hessian term in the chain rule. Choose the ratio interval so that both differential inequalities hold
wherever $H(B_c,B_e)$ differs from $\min(B_c,B_e)$. We obtain
\begin{equation*}
  B_*\asymp r^{-b},\qquad J_{\Gamma_*}B_*\le-cr^{-b-2},
\end{equation*}
equal to the selected exterior barrier outside a fixed compact set.
Only the compact problem \eqref{glu:corebarrierproblem} was used.

Transfer $B_*$ first to $\Sigma_{1,R}$. On the transition annulus the
scaled geometric difference tends to zero with $R$ by
\eqref{glu:fastcore}, so the strict differential inequality persists;
on the unchanged core it holds exactly. For fixed $R$, use the normal identification to define this
function on $\Sigma_{\sigma,R}$ and join it to the fractional
exterior barrier. The weighted compact convergence requirements for its
first two derivatives and for the geometry are included in the diagonal
choice of Lemma~\ref{glu:diagonal}. They preserve the strict negative
margin on the compact part and on the joining annulus.

The two barriers are normalized to have the same leading coefficient
on the end. Thus their cutoff commutator has a vanishing relative
coefficient on the classical transition annulus; on each fixed annulus
its fractional counterpart is small by compact convergence. No growing
multiplicative constant is introduced when the gluing radius increases.

Before applying \eqref{glu:concavemin}, define both functions $B_1,B_2$
on the whole surface. Choose positive extensions with bounds $Cr^{-b}$
and the required scaled derivatives, so that $B_2>2B_1$ on the interior region away from the collar, and
$B_1>2B_2$ on the exterior region away from it. Then $H$ differs from the minimum only
where both functions satisfy their strict differential inequalities.
The operator $LB_i(X)$ includes values $B_i(Y)$ from the whole
surface, including points where $H(B_1(Y),B_2(Y))$ equals the other
function. The additional integrals are bounded by
\eqref{glu:innerbarriertail} and its outer counterpart, and have size
$C\eps r^{-b-\tau}$. They are absorbed in the strict inequality.
The extensions and constants are fixed before choosing the diagonal
radius, so the final constants are independent of that radius.
For $b=\beta$
we obtain
\begin{equation}\label{glu:globalstrict}
  cr^{-\beta}\le B_\sigma\le Cr^{-\beta},\qquad
  L_\sigma B_\sigma\le-cr^{-\beta-\tau}.
\end{equation}
Repeating the construction with $b_s$ gives a slower strict barrier
$B_{s,\sigma}\asymp r^{-b_s}$.

\subsection{The weighted inverse}\label{glu:weighted-inverse}

In \cite[Section 7, (7.2) and Lemma 7.1]{DDPW}, reflection removes
the bounded odd translation field $Z_1$ of the classical catenoid.
The even logarithmic dilation field $Z_2$ enters the projected decaying
problem (7.9), coupled to the exterior problem through
(7.11)--(7.12). Here the exact classical graph has the positive
vertical field $w_*$ used in \eqref{glu:corebarrierproblem}; its slow
decay is excluded by weight $\beta$. On the nonminimal fractional
approximation the vertical normal is not a Jacobi field. The strict
barrier $B_\sigma$ replaces the exact positive field in
\cite[Corollary B.1]{DDPW}, giving direct comparison through
\eqref{glu:picone}. On a nonminimal background the perimeter second
variation also contains the curvature term of
\cite[Proposition B.2, (B.2)]{DDPW}; our coercivity is established
directly for $L_\sigma$ by the strict barrier.

All hypotheses of Proposition~\ref{mod:leafinverse} are now verified.
The kernel is positive and symmetric, its potential satisfies
\eqref{glu:potentialbound}, and the strict barriers are
\eqref{glu:globalstrict} and $B_{s,\sigma}$. The weighted tails are
\eqref{glu:innerbarriertail} and its outer counterpart. Thus the
same variational exhaustion and ratio comparison give
\begin{equation*}
  \sup r^\beta|u|\le C\|f\|_{Y_{\sigma,\beta}},\qquad L_\sigma u=f.
\end{equation*}
The local estimate \eqref{glu:schauder} applies because these same
sums control the weighted complementary integrals after chart scaling,
including their center differences. Consequently
\begin{equation}\label{glu:inverseestimate}
  T_\sigma=L_\sigma^{-1}:Y_{\sigma,\beta}\longrightarrow X_{\sigma,\beta},
  \qquad\|T_\sigma f\|_{X_{\sigma,\beta}}\le C\|f\|_{Y_{\sigma,\beta}},
\end{equation}
with constants independent of $\sigma$ and the diagonal radius.
Uniqueness and preservation of symmetries follow from the same
proposition. We use the variation identities of \cite{DDPW}; its catenoid inverse of Proposition 2.2 has different geometry
and weights. Our barriers and coarea estimates provide the required
verification here.

\subsection{The nonlinear correction and the graph property}\label{glu:nonlinearsection}

We first check that every small normal deformation in the fixed-point
ball defines an entire graph. Use the
horizontal coordinate $x$ of the background graph, and write
$w=(\nu_{\Sigma_\sigma})_9$. Since the norm of the differential of
$x\mapsto(x,F_{\Sigma_\sigma}(x))$ is $w^{-1}$, intrinsic derivatives
give
\begin{equation}\label{glu:projectionderivative}
  |D_x(h\nu_{\rm hor})|
  \le Cw^{-1}(|\nabla h|+|h||A|)
  \le C\|h\|_{X_{\sigma,\beta}}r^{m_\sigma-\beta-1}.
\end{equation}
The exponent is negative uniformly in the parameter. Choose
$\delta_0>0$ so small that this derivative is less than $1/2$ whenever
$\|h\|_{X_{\sigma,\beta}}\le\delta_0$. For every prescribed
$y\in\R^8$, the equation
\[
x=y-h(x)\nu_{\rm hor}(x)
\]
is a contraction on the complete space $\R^8$. Thus the horizontal
projection $\Pi_h(x)=x+h(x)\nu_{\rm hor}(x)$ is one-to-one and onto,
with uniformly bounded inverse differential. Its displacement is at
most $C\delta_0r^{-\beta}$, so the old and new horizontal weights are
comparable. The full normal map is an embedding. Its intrinsic
Jacobian stays between fixed positive bounds, since
$|h||A|+|\nabla h|\le C\delta_0r^{-\beta-1}$.

We evaluate the curvature at $X+h(X)\nu(X)$ and write it as a function
on $\Sigma_\sigma$:
\[
\mathcal H_\sigma(\Sigma_\sigma[h])
=E_\sigma+L_\sigma h+Q_\sigma(h),\qquad
E_\sigma=\mathcal H_\sigma(\Sigma_\sigma).
\]
Here $Q_\sigma$ is the explicit integral \eqref{eq:explicit-Q} with
$M=\Sigma_\sigma$ and $u=h$, evaluated at the displaced point.
We take the upward normal. In this
subsection abbreviate $X=X_{\sigma,\beta}$ and $Y=Y_{\sigma,\beta}$.
The expansion is exactly the normalized form of
\cite[(2.9), Proposition B.1]{DDPW}, as recorded in
\eqref{eq:DDPW-conventions}. Its local estimate uses their Section 8
calculation; the global estimate requires a different geometric bound.
In the catenoid problem the two conical sheets have separation
$O(\sqrt\eps\,r)$, and the normal cutoff in \cite[(8.2)]{DDPW} has
that thickness. Its normal derivatives cost
$(\sqrt\eps\,r)^{-1}$. The outer estimate in their Lemma 8.2 has
constant $C\eps^{-3/2}$ before normalization, yielding the
$\eps^{-1/2}$ loss in \cite[Proposition 2.3, (2.15)]{DDPW}.
Here Proposition~\ref{glu:core} and Lemma~\ref{glu:diagonal} give
injective normal coordinates for $|s|<cr$, uniformly in $\sigma$.
We can localize at this fixed fraction of $r$. The far interactions
are then controlled by the weighted area moments below; a
parameter-dependent thin cutoff is unnecessary.

\begin{proposition}[Uniform quadratic estimate]\label{glu:quadratic}
  There are constants $\delta_0,C>0$, independent of the diagonal radius
  and of $\sigma$ close to one, such that
  \begin{equation}\label{glu:qlipschitz}
    \|Q_\sigma(h)-Q_\sigma(k)\|_Y
    \le C(\|h\|_X+\|k\|_X)\|h-k\|_X
  \end{equation}
  whenever $\|h\|_X+\|k\|_X\le\delta_0$.
\end{proposition}

\begin{proof}
  At any intermediate graph, consider the second height variation in
  directions $v_1,v_2$. The local estimate of
  Section~\ref{sec:preliminaries}, obtained from
  \cite[Lemmas 8.3--8.5]{DDPW} with the exponents specified in
  Section~\ref{sec:preliminaries}, applies
  on the uniform tangent charts. The rescaled heights have norm
  $\|v_i\|_Xr^{-\beta-1}$, and curvature scales by $r^{-\sigma}$.
  Their bilinear contribution therefore has weight
  $r^{-2\beta-\tau-1}$. Taylor's integral formula gives
  \begin{equation}\label{glu:localquad}
    \|Q_\sigma(h)-Q_\sigma(k)\|_{C^\alpha_{\rm sc};\,\mathrm{loc}}
    \le C(\|h\|_X+\|k\|_X)\|h-k\|_Xr^{-2\beta-\tau-1}.
  \end{equation}
  For the remote part use Lemma~\ref{res:coarea}, transferred through
  the core insertion and the small normal maps:
  \begin{equation}\label{glu:weightedarea}
    \int_{\Sigma_\sigma\cap B_D(X)}r(Y)^{-b}\dd S_Y\le CD^{8-b},
    \qquad D\ge r(X)/2,\quad 0<b<7.
  \end{equation}
  The fixed compact core is harmless, the end insertion has a uniformly
  small normal derivative, and the area Jacobians and radii of the small
  graphs are comparable by \eqref{glu:projectionderivative}.
  In normal coordinates for the region between $M$ and $M[u]$, its
  volume Jacobian satisfies $|\partial_s j(Y,s)|\le C/r(Y)$.
  Thus two source variations give $|v_1(Y)v_2(Y)|/r(Y)$, while a
  kernel derivative at the source point or the evaluation point costs
  one inverse separation.
  Outside the tangent ball the resulting majorants are
  \begin{align*}
    &\eps\int
    \frac{|v_1(Y)v_2(Y)|/r(Y)}{|X-Y|^{8+\tau}}\dd S_Y,\\
    &\eps\int
    \frac{|v_1(X)v_2(Y)|+|v_2(X)v_1(Y)|+|v_1(Y)v_2(Y)|}
    {|X-Y|^{9+\tau}}\dd S_Y,\\
    &|v_1(X)v_2(X)|r^{-\tau-1}.
  \end{align*}
  For the last term use $|\chi_{E^c}-\chi_E|\le1$:
  \[
  C\eps |v_1(X)v_2(X)|
  \int_{|Z-X|>cr}|Z-X|^{-11-\sigma}\dd Z
  \le C\eps |v_1(X)v_2(X)|r^{-2-\sigma}.
  \]
  The localizing cutoff has derivatives of scale $r^{-1}$ and the
  same majorants. The small normal displacements preserve every remote
  separation and tube Jacobian. Comparable annuli between distances
  $cr$ and $r/2$ give the asserted power directly.
  For the remaining annuli suppress the factor
  $\|v_1\|_X\|v_2\|_X$ and use ambient dyadic annuli of radii
  $D_j=2^jr/2$, $j\ge0$. For the first integral,
  \eqref{glu:weightedarea} with $b=2\beta+1$ gives
  \[
  C\eps\sum_{j\ge0}D_j^{-8-\tau}D_j^{8-(2\beta+1)}
  \le C\eps r^{-2\beta-\tau-1}.
  \]
  The second integral's source-source term uses $b=2\beta$.
  For a mixed term the remaining sum is
  \[
  C\eps r^{-\beta}\sum_{j\ge0}
  D_j^{-9-\tau}D_j^{8-\beta}
  \le C\eps r^{-2\beta-\tau-1}.
  \]
  The largest required source weight is $2\beta+1=6<7$.
  Equivalently, source annuli $\rho<r/2$ in the first integral give,
  after vertical integration, the moment
  \[
  C\eps r^{-8-\sigma}\sum_{\rho<r/2}\rho^{6-2\beta}
  \le C\eps r^{-2\beta-\tau-1}.
  \]
  This is the restriction on the smaller source radii; the ambient
  annular sums also control large distances.

  For center differences, a remote kernel derivative on the $j$th
  annulus costs $Cr/D_j\le C$ in scaled coordinates. The cutoff and
  evaluation heights obey the same scaled H\"older bounds. The absolute
  moments therefore give the remote $C^\alpha$ estimate with the same
  power. These bounds justify differentiation first with truncated
  kernels, and then removal of both the diagonal and ambient cutoffs.
  Apply the integral Taylor formula as in
  \cite[proof of Lemma 8.1]{DDPW}, retaining the moving evaluation
  point throughout. The resulting full estimate is
  \eqref{glu:localquad}; multiplication by $r^{\beta+\tau}$ leaves
  $r^{-\beta-1}\le1$ and proves \eqref{glu:qlipschitz}.
\end{proof}

Set $e_\sigma=\|E_\sigma\|_Y\to0$, as in
Lemma~\ref{glu:diagonal}. Let $C_T$ be the inverse bound in
\eqref{glu:inverseestimate} and take the ball of radius
$\rho_\sigma=2C_Te_\sigma$ in $X$. For
\begin{equation}\label{glu:fixedpoint}
  \mathcal T_\sigma(h)=-T_\sigma(E_\sigma+Q_\sigma(h))
\end{equation}
we have
$\|\mathcal T_\sigma(h)\|_X\le C_Te_\sigma+C\rho_\sigma^2$
and Lipschitz constant at most $2C\rho_\sigma$. For $\sigma$ close to
one, $\rho_\sigma<\delta_0$ and the map is a contraction of this ball
into itself. Its fixed point satisfies
\begin{equation}\label{glu:smallcorrection}
  \|h_\sigma\|_X\le2C_Te_\sigma\longrightarrow0.
\end{equation}
The corrected hypersurface has zero fractional curvature.

\begin{proposition}[Preservation of the entire graph]\label{glu:graphpreservation}
  For $\sigma$ sufficiently close to one,
  $\Sigma_\sigma[h_\sigma]$ is a smooth entire graph over $\R^8$.
  It has the stated rotational and interchange symmetries and is not affine.
\end{proposition}

\begin{proof}
  The horizontal inverse was proved above for every height of norm
  at most $\delta_0$. The vertical normal remains positive:
  writing $\delta=\|h_\sigma\|_X$, the normal-graph formula gives
  \[
  (\nu_{\Sigma_\sigma[h_\sigma]})_9
  \ge cr^{-m_\sigma}-C\delta r^{-\beta-1}
  \ge\tfrac c2r^{-m_\sigma}>0.
  \]
  Here $|h_\sigma||A|$ is small and
  $|\nabla h_\sigma|\le C\delta r^{-\beta-1}$, while
  $m_\sigma<\beta+1$ uniformly.

  At a fixed new horizontal coordinate the vertical change has size
  \begin{equation}\label{glu:verticalcorrection}
    |F_{\Sigma_\sigma[h_\sigma]}-F_{\Sigma_\sigma}|
    \le C\delta r^{m_\sigma-\beta}.
  \end{equation}
  Use the fixed-horizontal-point velocity identity from the proof of
  \eqref{glu:verticalnormalvariation}, now along $s h_\sigma$.
  The normal velocity is bounded by $C\delta r^{-\beta}$, and the
  preceding normal estimate gives $w_s\asymp r^{-m_\sigma}$ throughout
  this path. The inverse projections preserve comparable radii.
  Integration gives \eqref{glu:verticalcorrection}; its decay preserves
  the superlinear principal growth.

  The equation, the approximation and its inverse are equivariant under
  the prescribed rotations and the interchange $(u,v)\mapsto(v,u)$. Applying any such
  symmetry to the small solution gives another fixed point in the same
  ball, with the corresponding signed normal height. Uniqueness preserves
  the symmetries. For each fixed $\sigma$, spatial difference quotients
  of the graph equation satisfy uniformly elliptic equations of order
  $\tau$ on rotated local charts. The local Schauder estimate and
  iteration give smoothness, as in the graph regularity argument of
  \cite{BFV}. The nonzero superlinear profile identified below proves
  that this entire graph is not affine.
\end{proof}

\begin{proof}[Proof of Theorem~\ref{thm:main}]
  The exact fractional leaf and the cylindrical matching give the
  principal end. Proposition~\ref{glu:core} and
  Lemma~\ref{glu:diagonal} insert the compatible classical core with
  $e_\sigma\to0$. The inverse \eqref{glu:inverseestimate}, the quadratic
  estimate \eqref{glu:qlipschitz} and the fixed point
  \eqref{glu:fixedpoint} yield the correction
  \eqref{glu:smallcorrection}. Its normalized and physical curvatures
  vanish simultaneously, and Proposition~\ref{glu:graphpreservation}
  gives a smooth entire graph with the required symmetries.

  Let $\Lambda_\sigma(x)$ be the signed scale of the exact leaf through
  $x$. The foliation by dilates makes this a continuous function off the
  cone, homogeneous of degree one. Define
  \[
  F_{\infty,\sigma}(x)
  =a_0\operatorname{sgn}(\Lambda_\sigma(x))
  |\Lambda_\sigma(x)|^{p_\sigma}.
  \]
  It is homogeneous of degree $p_\sigma=1+m_\sigma$, and nonzero off
  $C$. As a point on a fixed horizontal sphere approaches $C$, the
  leading conical expansion of the leaf gives
  $\operatorname{dist}(x,C)\asymp|\Lambda_\sigma(x)|^{p_\sigma}$.
  Thus this profile extends continuously by zero to $C$.

  Fix a compact set $K\subset\R^8\setminus C$. The absolute value of
  this profile is bounded below on $K$. Define $t_L(x)$ by
  $A_{\rm out}(t_L(x))=L^{p_\sigma}F_{\infty,\sigma}(x)$.
  The signed scale of the corresponding outer leaf is
  $L\Lambda_\sigma(x)$, so its horizontal slice contains $Lx$.
  Moreover $|t_L|\asymp L^{p_\sigma}$ uniformly on $K$.
  Since $p_\sigma>a$, these candidate points lie in the unchanged
  outer region for large $L$ and hence give its actual graph values.
  Thus its defining label identity is
  \[
  A_{\rm out}(F_{\mathcal P_\sigma}(Lx))
  =L^{p_\sigma}F_{\infty,\sigma}(x).
  \]
  Write $\theta_\sigma=q-d<1$. Because
  $A_{\rm out}(t)=t+O(|t|^{\theta_\sigma})$ and
  $A_{\rm out}'(t)\to1$, inversion of this identity gives, uniformly on
  $K$,
  \[
  L^{-p_\sigma}F_{\mathcal P_\sigma}(Lx)
  =F_{\infty,\sigma}(x)
  +O(L^{-p_\sigma(1-\theta_\sigma)}).
  \]
  The compact core insertion is absent there, and
  \eqref{glu:verticalcorrection} contributes at most
  $C\delta L^{m_\sigma-\beta-p_\sigma}
  =C\delta L^{-\beta-1}$ after this division. This proves
  \eqref{eq:mainasymptotic}; the degree tends to three as
  $\sigma\uparrow1$.

  Finally, Lemma~\ref{glu:diagonal} gives
  $\Sigma_\sigma\to\Gamma_*$ in $C^2$ on every fixed compact set.
  On such a set the norm $X_{\sigma,\beta}$ controls the
  $C^{2,\alpha_\sigma}$ norm of $h_\sigma$, with
  $\alpha_\sigma\ge\alpha/2$. Hence \eqref{glu:smallcorrection} gives
  convergence of the normal maps in $C^2$. The horizontal projections are $C^2$ maps converging to the identity:
  the terms in their first two derivatives involve only the first two
  derivatives of $h_\sigma$ and the fixed smooth geometry on that
  compact set. Their first derivatives are uniformly invertible by
  \eqref{glu:projectionderivative}, so the inverse maps converge in
  $C^2$ there as well. The entire graph functions therefore
  converge locally in $C^2$ to the selected classical graph $F_*$.
\end{proof}

\section*{Acknowledgement}

The research of M. del Pino is supported by the Royal Society Research Professorship grant RP-R1-180114 and by the ERC/UKRI Horizon Europe grant ASYMEVOL, EP/Z000394/1.  The research of J. Wei is partially supported by GRF of RGC of Hong Kong entitled ``On critical and supercritical Fujita equation''.  This project originated in 2013 and has been periodically updated since then, with the core ideas and technical framework established in an earlier manuscript. During the final preparation of this manuscript, the authors used ChatGPT (OpenAI) to assist with language editing and exposition. All mathematical arguments, proofs, and conclusions were independently verified by the authors, who take full responsibility for the content of this paper.

\end{document}